\documentclass[11pt]{amsart}
\usepackage{wasysym}
\usepackage{xypic}
\usepackage{enumitem}
\usepackage{fontawesome}
\usepackage[hidelinks]{hyperref}

\usepackage{xr}
\newfont{\wcyr}{wncyr10 at 10pt}
\newfont{\wcyi}{wncyi10 at 10pt}
\newfont{\swcyi}{wncyi10 at 8pt}
\newfont{\wcyb}{wncyb10 at 10pt}

\newenvironment{Corollary}%
   {\smallskip \noindent {\sc Corollary.}}%
   {\bigskip}
\newenvironment{Definition}%
   {\smallskip \noindent {\sc Definition.}}%
   {}
\newenvironment{Example}%
   {\smallskip \noindent {\sc Example.}}%
   {\bigskip}
\newenvironment{Examples}%
   {\smallskip \noindent {\sc Examples.}}%
   {\bigskip}
\newenvironment{Proof}%
   {\noindent {\sc Proof:}}%
   {\hfill\raisebox{-2pt}{\faThumbsOUp} \bigskip}
\newenvironment{Proposition}%
   {\smallskip \noindent {\sc Proposition.}}%
   {\bigskip}
\newenvironment{Theorem}%
   {\smallskip \noindent {\sc Theorem.}}%
   {\bigskip}
\newenvironment{Lemma}%
   {\smallskip \noindent {\sc Lemma.}}%
   {\bigskip}
   {\smallskip \noindent {\sc Remark.}}%
   {\bigskip}
   {\smallskip \noindent {\sc Remarks.}}%
   {\bigskip}
\newenvironment{Question}%
   {\smallskip \noindent {\sc Question.}}%
   {\bigskip}
   {\smallskip \noindent {\sc Questions.}}%
   {\bigskip}
   {\smallskip \noindent {\sc Conjecture.}}%
   {\bigskip}

\newcommand{\cB}{{\mathcal B}}
\newcommand{\bbC}{\mathbb{C}}

\newcommand{\cG}{{\mathcal G}}
\newcommand{\calH}{{\mathcal H}}

\newcommand{\cK}{{\mathcal K}}

\newcommand{\N}{\mathbb{N}}
\newcommand{\Q}{\mathbb{Q}}
\newcommand{\R}{\mathbb{R}}

\newcommand{\cN}{{\mathcal N}}

\newcommand{\cP}{{\mathcal P}}

\newcommand{\bbP}{{\mathbb P}}
\newcommand{\cQ}{{\mathcal Q}}

\newcommand{\cS}{{\mathcal S}}

\newcommand{\cV}{{\mathcal V}}

\newcommand{\cX}{{\mathcal X}}

\newcommand{\cY}{{\mathcal Y}}

\newcommand{\cZ}{{\mathcal Z}}

\DeclareMathOperator{\diam}{diam}
\DeclareMathOperator{\supp}{supp}
\DeclareMathOperator{\Span}{span}

\newcommand{\Set}{{\mathbf {Set}}}

\DeclareMathOperator{\Fin}{Fin}

\newcommand{\imp}{\Rightarrow}
\newcommand{\ifff}{\Longleftrightarrow}

\newcommand{\cstar}{$\mathrm{C}^*$}

\DeclareMathOperator{\Add}{Add}

\newcommand{\bbN}{\N}

\DeclareMathOperator{\dom}{dom}
\DeclareMathOperator{\Seq}{Seq}
\DeclareMathOperator{\Ran}{Ran}
\newcommand{\Hilbert}{\bf {Hilbert}}
\newcommand{\AWstar}{$\mathrm{AW}^*$}

\DeclareMathOperator{\Sym}{Sym}
\DeclareMathOperator{\HSym}{HSym}
\DeclareMathOperator{\id}{id}
\newenvironment{ProofOf}[1]%
{\noindent {\sc #1:}}%
{\hfill\raisebox{-2pt}{\faThumbsOUp} \bigskip}

\newenvironment{Problem}%
{\smallskip \noindent {\sc Problem.}}%
{\bigskip}

\title{Hilbert Spaces Without The Countable Axiom of Choice}
\author{Bruce Blackadar, Ilijas Farah, and Asaf Karagila}

\thanks{I.F. was partially supported by NSERC}
\address[B.B]{Department of Mathematics and Statistics,
University of Nevada, Reno\\
Reno, NV, USA  89557}
\email{bruceb@unr.edu}

\address[I.F]{Department of Mathematics and Statistics,
	York University,
	4700 Keele Street,
	Toronto, Ontario, Canada, M3J
	1P3\\
Matemati\v cki Institut SANU\\
	Kneza Mihaila 36\\
	11\,000 Beograd, p.p. 367\\
	Serbia}
\email{ifarah@yorku.ca}
\address[A.K.]{School of Mathematics\\ University of Leeds\\ Leeds, LS2 9JT, UK}
\email{karagila@math.huji.ac.il}
\thanks{A.K. was supported by a UKRI Future Leaders Fellowship [MR/T021705/2]}
\thanks{We would like to thank the anonymous referee for a very useful report.}

\date{\today}

\begin{document}

\maketitle

\begin{center}
{\em To the memory of Eberhard Kirchberg}
\end{center}

This article examines Hilbert spaces constructed from sets whose existence is incompatible
with the Countable Axiom of Choice (CC).  Our point of view is twofold:

(1)  We examine what can and cannot be said about Hilbert spaces and operators on them
in ZF set theory without any assumptions of Choice axioms, even the CC.  It turns out
that quite a bit can be said in this context, some of which seems a little counterintuitive.  We are
especially interested in Hilbert spaces whose dimension is ``finite'' in a sense incompatible
with the CC, where strange properties sometimes hold.  One philosophical consequence
of not assuming CC is that Hilbert spaces come in a great many more varieties and
are not nearly as homogeneous as normally thought (and are much more interesting!)

(2)  We view Hilbert spaces as ``quantized'' sets and obtain some set-theoretic results from
associated Hilbert spaces.

The article is written in a partly expository style.  Some of the results were previously obtained,
primarily by N.\ Brunner and coauthors in \cite{BrunnerLinear}, \cite{BrunnerHilbertraume}, and \cite{BrunnerSB} (although these
papers have a somewhat different point of view and contain some errors; see the discussion in \S\ref{RussellOrthSeq}).
Many of our basic results are modest and qualify as ``low-hanging fruit''; we have certainly not
picked all the fruit and invite other interested readers to join the harvest.

A subsequent article \cite{BlackadarF} will discuss C*-algebras without AC, including the validity
of the two Gelfand-Naimark representation theorems.

%It is perhaps not too surprising that the behavior of (nonseparable) Hilbert spaces in the absence of the Axiom of Choice resembles that of arbitrary (separable) Banach spaces. For example, Proposition~\ref{Q(H)OneDimensional} gives a Hilbert space which is hereditarily indecomposable (Definition~\ref{HilbDFDef} (v),  cf. \cite{GowMa:Unconditional}) every operator on which is scalar+compact (cf. \cite{argyros2011hereditarily}).

{\bf Caution:}  Without
CC one must be very careful since many familiar techniques we use, sometimes without
much thought, do not work.  The principal tools we cannot use are Zorn's Lemma, Tikhonov's
Theorem, the Hahn-Banach Theorem, and the Baire Category Theorem, all of which require
some version of the AC.  (Note, however, that we can still make finitely many choices.)  Thus we will have to prove or reprove some ``standard'' or ``well-known''
properties of Hilbert spaces.  So we cannot assume we know anything about Hilbert spaces beyond
some basic facts obviously independent of Choice, such as the CBS inequality and
Parallelogram Law.  \cite{BlackadarReal} is a convenient reference for some standard arguments,
which we will use as standard reference, but they can also be found in many other references.
We use the capitalized word {\em Choice} to generically refer to versions of the Axiom of Choice
(e.g.\ ``can be proved without any Choice'' means ``theorem of ZF''),
and we will not assume any Choice without explicit mention.

\tableofcontents

\subsection*{Brief Summary} 
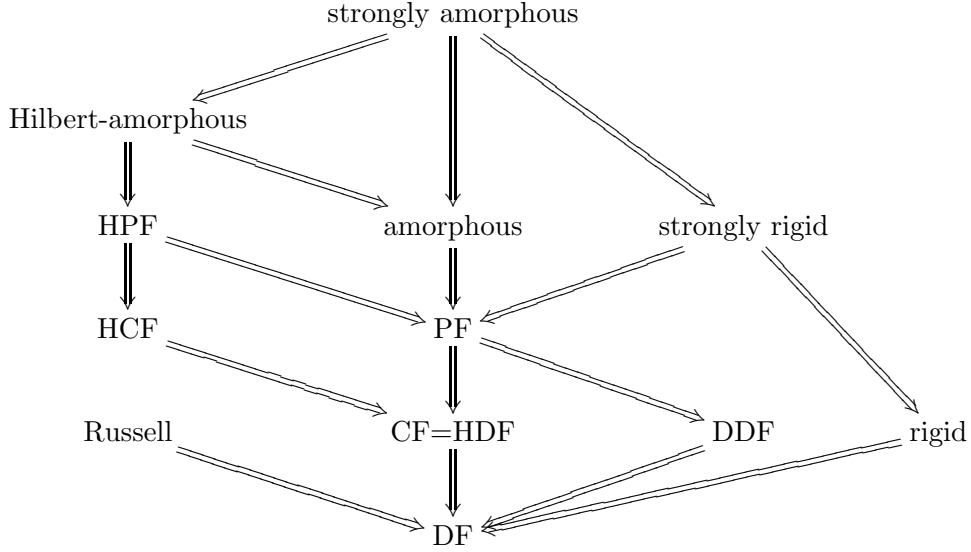
\begin{figure}[h]

	$$\xymatrix{ & \mbox{strongly amorphous}\ar@{=>}[dl]\ar@{=>}[dd]\ar@{=>}[ddr] & & \\
		\mbox{Hilbert-amorphous}\ar@{=>}[d]\ar@{=>}[dr] & & & \\ \mbox{HPF}\ar@{=>}[d]\ar@{=>}[dr] & \mbox{amorphous}\ar@{=>}[d] & \mbox{strongly rigid}\ar@{=>}[ddr]\ar@{=>}[dl] & \\
		\mbox{HCF}\ar@{=>}[dr] & \mbox{PF}\ar@{=>}[d]\ar@{=>}[dr] & & \\ \mbox{Russell}\ar@{=>}[dr] & \mbox{CF=HDF}\ar@{=>}[d] & \mbox{DDF}\ar@{=>}[dl] & \mbox{rigid}\ar@{=>}[dll] \\
		& \mbox{DF} & & }$$
	\caption{\label{Fig.DF2}Implications between different flavors of Dedekind--finite sets discussed in this section. Hilbert-amorphous, HPF, HCF, and HDF will be defined in \S\ref{S.HDF}.}
\end{figure}

In Section \ref{S.Completeness} we introduce the notion of $\sigma$-completeness of metric spaces and argue that it is the correct notion of completeness in a context where CC fails (see also Example~\ref{CauchyCompleteEx}). In Secion \ref{S.Hilbert} we establish the theory of Hilbert spaces in this context. Spaces $\ell^2(X)$ are introduced in Section~\ref{S.Orthonormal}, and in Section~\ref{S.FiniteSeparable} we show that the analysis of finite-dimensional and separable Hilbert spaces is unaffected even if CC fails. Dedekind-finite sets and their variations are introduced in Section~\ref{S.DF+}. The most important notion for us is that of a Cohen-finite set (Definition~\ref{Def.CF}), and a CF set is a Cohen-finite, infinite set.  Simple permanence properties of CF sets are established in Subsection~\ref{S.Permanence}. This includes Lemma~\ref{L.BackAndForth}, used throughout. Fig.~\ref{Fig.DF} represents implications between different flavors of Dedekind-finite sets. In Section~\ref{S.ell2} we establish basic results on $\ell^2(X)$ spaces for a DF (Dedekind-finite and infinite) set $X$. In Proposition~\ref{StrongDFFiniteSupport} we prove, among other equivalences, that $X$ is Cohen-finite if and only if the standard basis is a Hamel basis. In Theorem~\ref{L.ExistenceOfBasis} we find (using a rather peculiar variant of a Dedekind-finite set, proved consistent with ZF by J. Truss) a Hilbert space without an orthonormal basis, and even without an infinite orthonormal sequence, and in Example~\ref{NonuniqueBasis} we prove that if $X$ is a DF set which is not Cohen-finite then $\ell^2(X)$ has orthonormal bases of different cardinalities. The question whether every Hilbert space is (provably in ZF) isometric to a subspace of one with an orthonormal basis is discussed in Subsection~\ref{S.universality}.
In Section~\ref{S.Bounded}  we turn to the study of bounded linear operators, and in Section~\ref{S.B(H)} we consider the algebra of bounded linear operators on $\ell^2(X)$ for a DF set $X$. If $X$ is CF then   this algebra is non-separable, but it is stably finite (Corollary~\ref{C.XCFstablyfinite}) and every compact operator has finite rank (Proposition~\ref{CompactFinRankProp}). In Section~\ref{S.Q(H)} we study the Calkin algebra $\cQ(\ell^2(X))$ and prove that it is one-dimensional if $X$ is strongly amorphous (Proposition~\ref{Q(H)OneDimensional}) and that it can be non-separable and abelian (Proposition~\ref{P.Calkin.abelian}) or nonabelian and stably finite (Proposition~\ref{S.Q.nonabelian}) for specific choices of a CF set $X$. In Corollary~\ref{C.Baire} we prove that in Solovay's model there is a stably finite \cstar-algebra with no tracial states. This example is commutative and it has no nonzero representation  on a Hilbert space. In  Proposition~\ref{P.stably-finite-no-tracial-states} we show that if there is a CF set~$X$ such that its power set has no finitely additive probability measure that vanishes on singletons, then $\cB(\ell^2(X))$ is stably finite and has no tracial states, and is therefore not satisfactory (this is an immediate  consequence of a well-known fact that if all sets of reals have the Property of Baire then there is no probability measure on $\cP(\bbN)$ that vanishes on the singletons).
In Section~\ref{S.Spectrum} we prove that for a power Dedekind-finite set~$X$, the spectrum of every bounded linear operator $T$ on $\ell^2(X)$ is finite and every element of the spectrum is an eigenvalue. In Section~\ref{S.HDF} we introduce Hilbert space analogs of principal flavors of Dedekind-finite sets (see Fig.~\ref{Fig.DF2}).  Models of ZF used in the previous sections are constructed in Section~\ref{S.ModelsofZF}, and we conclude with a short list of open problems in Section~\ref{S.Problems}.

\section{Completeness and its Consequences}\setcounter{paragraph}{0}\label{S.Completeness}

%In \cite{}, an example in the absence of the Countable AC is given of a nonreflexive Hilbert space.  While
%this example is very interesting, it is only possible because the ``wrong'' definition of completeness is used:
The first technicality concerns the ``right'' notion of completeness for Hilbert spaces and, more
generally, metric spaces.
The Cauchy sequence definition (which we will call {\em Cauchy completeness}) turns out not to be useful in the absence of CC.
Instead, the following version of completeness should be used.

\paragraph{}
\begin{Definition}\label{ACMetricCompDef}
Let $(X,\rho)$ be a metric space.

\smallskip

\noindent
(i)  $(X,\rho)$ is {\em absolutely closed} if, whenever $\phi$ is an isometric embedding of $(X,\rho)$
into a metric space $(Y,\rho')$, $\phi(X)$ is closed in $(Y,\rho')$.

\smallskip

\noindent
(ii)  $(X,\rho)$ is {\em $\sigma$-complete} if, whenever $(A_n)$ is a decreasing sequence of
nonempty closed  subsets of $X$ with $\diam(A_n)\to0$, then $\cap_n A_n$ is nonempty
(it is then necessarily a singleton). In \cite{BrunnerSB} spaces with this property were called \emph{Cantor complete}.

\smallskip

\noindent
(iii)  $(X,\rho)$ is {\em totally complete} if, whenever $(A_i)$ is an indexed set of nonempty closed subsets of $X$, directed by reverse inclusion (i.e.\ $A_j\subseteq A_i$ if and only if
$i\leq j$, and for every $i,j$ there is a $k$ with $A_k\subseteq A_i\cap A_j$; the indices run over an arbitrary directed set) with $\diam(A_i)\to0$
(i.e.\ for every $\epsilon>0$ there is an $i_0$ with $\diam(A_{i_0})<\epsilon$ and hence
$\diam(A_i)<\epsilon$ for all $i>i_0$),
then $\cap_i A_i$ is nonempty
(it is then necessarily a singleton).

\smallskip

\noindent
(iv)  $(X,\rho)$ is {\em uniformly complete} if the uniform space defined by the metric is
complete,  in the sense that every Cauchy net is convergent.
\end{Definition}

\bigskip

In the context of the Axiom of Choice, a metric space is uniformly complete if and only if it is Cauchy complete. However, without CC in general a non-convergent Cauchy sequence cannot be extracted from a non-convergent Cauchy net (cf.\ \ref{StrongDFFiniteSupport}, \ref{L.SDFBasis}).
The following is a theorem of ZF, i.e.\ no form of Choice is used.  If the CC is assumed, all five conditions are equivalent.  See \cite[XII.16.9.5]{BlackadarReal} for a proof.

\paragraph{}
\begin{Theorem}\label{ACCompleteThm}
Let $(X,\rho)$ be a metric space.  The following are equivalent:
\begin{enumerate}
\item[(i)]  $(X,\rho)$ is absolutely closed.
\item[(ii)]  $(X,\rho)$ is $\sigma$-complete.
\item[(iii)]  $(X,\rho)$ is totally complete.
\item[(iv)]  $(X,\rho)$ is uniformly complete.
\end{enumerate}
These conditions imply
\begin{enumerate}
\item[(v)] $(X,\rho)$ is Cauchy complete.
\end{enumerate}
\end{Theorem}

\smallskip

We will use the term $\sigma$-complete generically for a metric space satisfying (i)--(iv).  A closed
subset of a $\sigma$-complete metric space is obviously $\sigma$-complete.  A compact metric
space is $\sigma$-complete.

\smallskip

Since every uniform space has an essentially unique completion (cf.\ e.g.\ \cite[\S XII.16.8]{BlackadarReal},
and note that no Choice is used there), we obtain:

\paragraph{}
\begin{Corollary}\label{MetCompletion}
If $(X,\rho)$ is a metric space, there is a canonical $\sigma$-complete metric space $(\bar X, \bar\rho)$
containing $(X,\rho)$.  If $\phi:X\to Y$ is an isometry from $X$ to a $\sigma$-complete metric
space $(Y,\tau)$, then $\phi$ extends uniquely to an isometry from $\bar X$ to the closure of
$\phi(X)$ in $Y$.  Thus the completion of $(X,\rho)$ is unique up to isometry which is the identity
on~$X$.
\end{Corollary}

\section{Hilbert Spaces}\label{S.Hilbert}

\paragraph{}
\begin{Definition}\label{HilbSpaceDef}
A {\em Hilbert space} is a complex inner product space which is $\sigma$-complete under the induced
norm (almost everything we do works at least as well for real Hilbert spaces).  More generally, a {\em (real or complex) Banach space} is a $\sigma$-complete normed
(real or complex) vector space.
\end{Definition}

\smallskip

In the presence of CC, this agrees with the usual definition, and we believe it is the ``right''
definition in the general case.  Note that a Hilbert space by this definition is Cauchy-complete.
Every (complex) inner product space has a canonical $\sigma$-completion (\ref{MetCompletion}).  Since closed bounded
subsets of a finite-dimensional normed vector space are compact (no Choice is needed to prove
this), any finite-dimensional inner product space is a Hilbert space.  The problem with the
example in \cite{BrunnerSB} (cf.\ \S\ref{CauchyCompleteEx}) is that, while it is Cauchy complete, it is not $\sigma$-complete (and thus not a
Hilbert space by our definition).

\paragraph{}
The CBS inequality and the Parallelogram Law hold in any inner product space (no Choice
required), and the following fact (in ZF) is a fundamental consequence.  See e.g.\
\cite[XVI.9.4.2]{BlackadarReal} for a proof.

\paragraph{}
\begin{Lemma}\label{HilbCloseApproxLem}
Let $C$ be a nonempty convex set in an inner product space $\cV$, $\xi\in\cV$, and $\epsilon>0$.

If $\eta,\zeta\in C$ with\footnote{We write $\rho(\xi,C)=\inf\{\|\xi-\eta\|\mid \eta\in C\}$.}
	 $\|\xi-\eta\|^2<\rho(\xi,C)^2+\epsilon$ and $\|\xi-\zeta\|^2<\rho(\xi,C)^2+\epsilon$,
then $\|\eta-\zeta\|^2<4\epsilon$.
\end{Lemma}

\bigskip

The most important consequence of $\sigma$-completeness in inner product spaces is that the
following fundamental fact holds (in ZF!).  The usual proof uses Cauchy sequences and requires
CC, but there is an alternate proof from $\sigma$-completeness not requiring CC.  This proof
can be found in some references, but we give it here to illustrate the
use of $\sigma$-completeness.

\paragraph{}\label{CloseVectThm}
\begin{Theorem} {\sc [Closest Vector Property]}\label{HilbCloseApproxThm}
Let $ \calH$ be a Hilbert space and $C$ a nonempty closed convex subset of $ \calH$.  For any $\xi\in \calH$,
there is a unique $\eta\in C$ such that
$$\|\xi-\eta\|=\rho(\xi,\eta)=\rho(\xi,C)=\inf_{\zeta\in C}\|\xi-\zeta\|=\min_{\zeta\in C}\|\xi-\zeta\|\ .$$
\end{Theorem}

\begin{Proof}
Let $r=\rho(\xi,C)$.  For each $n\in\N$ let $A_n$ be the intersection of $C$ with the closed ball of
radius $r+\frac{1}{n}$ around $\xi$.  Then each $A_n$ is a nonempty closed convex set in $ \calH$, with
$A_{n+1}\subseteq A_n$, and $\rho(\xi,A_n)=\rho(\xi,C)$.  Lemma
\ref{HilbCloseApproxLem} implies that $\diam(A_n)\to 0$ as $n\to\infty$.  Thus by
$\sigma$-completeness $\cap_n A_n=\{\eta\}$ for some $\eta\in C$.
Clearly $\|\xi-\eta\|=r$ and $\eta$ is the unique vector in $C$ with this property.
\end{Proof}

A key consequence of \ref{CloseVectThm} is the existence of orthogonal complements:

\paragraph{}
\begin{Corollary}\label{OrthDecompProp}
Let $ \calH$ be a Hilbert space, and $\cY$ a closed subspace of $ \calH$.  Then every $\xi\in \calH$
can be uniquely written as $\eta+\zeta$, where $\eta\in\cY$ and $\zeta\in\cY^\perp$.  This $\eta$ is the
closest vector in $\cY$ to $\xi$.  So $\cY$ and $\cY^\perp$ are complementary closed
subspaces of $ \calH$, and in particular $\cY$ has an orthogonal complement, and $(\cY^\perp)^\perp=\cY$.  There is an orthogonal projection $P_{\cY}$ from $ \calH$ onto $\cY$ with null space $\cY^\perp$.
\end{Corollary}

%\begin{Proof}
%Existence of the closest vector $\eta$ comes from \ref{CloseVectThm}, and the usual proof
%\cite[XVI.9.5.7]{BlackadarReal} shows that $\zeta=\xi-\eta$ is orthogonal to $\cY$, and that the decomposition is unique.
%$P_{\cY}$ is defined by $P_{\cY}(\xi)=\eta$.
%\end{Proof}

\paragraph{}
\begin{Theorem} {\sc [Riesz Representation Theorem]}\label{RieszRepThm}
Let $ \calH$ be a Hilbert space, and $\phi$ a bounded linear functional on $ \calH$.  Then there is a
unique vector $\eta\in \calH$ with $\phi(\xi)=\langle \xi,\eta\rangle$ for all $\xi\in \calH$.  The map $\phi\mapsto \eta$
is a (conjugate-linear) isometry from $ \calH^*$ onto $ \calH$.  In particular, $ \calH$ is reflexive.
\end{Theorem}

\smallskip

The standard proofs work verbatim.

\bigskip

We also note the following immediate (actually somewhat nontrivial) consequence of \ref{MetCompletion}:

\paragraph{}
\begin{Proposition}\label{DenseIsom}
Let $ \calH$ and $ \calH'$ be Hilbert spaces, and $ \calH_0$ a dense subspace of $ \calH$.  An isometric
linear map $T$ from $ \calH_0$ to $ \calH'$ extends uniquely to an isometric linear map $\bar T$ from
$ \calH$ to $ \calH'$.  If $T( \calH_0)$ is dense in $ \calH'$, then $\bar T$ is surjective.
\end{Proposition}

\smallskip

More generally, since a bounded linear operator is uniformly continuous, any bounded
$T: \calH_0\to \calH'$ extends uniquely to $ \calH$.

\section{Orthonormal Bases and $\ell^2$-Spaces}\setcounter{paragraph}{0}\label{S.Orthonormal}

%\paragraph{}
One standard property of Hilbert spaces is dramatically absent in our non-Choice setting:
existence of orthonormal bases. Every Hilbert space that has a well-ordered dense subset has an orthonormal basis, even a well-ordered one (Theorem~\ref{L.l2(kappa)}). In particular, every separable Hilbert space has an orthonormal basis. A Hilbert space with no well-ordered dense subset is non-separable, and without the Axiom of Choice (or at least the axiom of Dependent Choices) it does not necessarily have a separable infinite-dimensional subspace (cf.\ \ref{StrongDFFiniteSupport}, \ref{L.SDFBasis}). Because of this, in the non-Choice setting a non-separable Hilbert space is not necessarily either  ``larger'' or ``smaller'' than a separable one).  Not every Hilbert space has an orthonormal basis without
some additional assumptions (\ref{L.ExistenceOfBasis}).
In this section, we will examine Hilbert spaces which come
equipped with a natural orthonormal basis, which might be better behaved.

\paragraph{}
The first thing we need to make precise is the notion of the sum of a function over a set.  We will
take it to mean the integral with respect to counting measure; thus the sum of a nonnegative
real-valued function on a set $X$ is the supremum of the sums over finite subsets of $X$.

\paragraph{}
\begin{Definition}\label{L2Def}
Let $X$ be a set.  Define
$$\ell^2(X)=\left \{ \eta:X\to\bbC:\sum_{x\in X}|\eta(x)|^2<\infty\right \} \ .$$
\end{Definition}

\paragraph{}
The usual proofs show that $\ell^2(X)$ is a complex vector space under pointwise operations,
and the CBS inequality shows that if $\eta,\zeta\in\ell^2(X)$, then the net
$$\langle \eta,\zeta\rangle_F=\sum_{x\in F}\eta(x)\overline{\zeta(x)}$$
in $\bbC$ indexed by finite subsets $F$ of $X$, directed by inclusion, converges to a complex number we
call $\langle \eta,\zeta\rangle$, and $\langle\cdot,\cdot\rangle$ is an inner product on $\ell^2(X)$.

\paragraph{}
\begin{Proposition}\label{L2CompleteProp}
Let $X$ be a set.  Then $\ell^2(X)$ is $\sigma$-complete, hence a Hilbert space.
\end{Proposition}

\begin{Proof}
The proof follows the pattern of the usual proof of completeness.  Let $(A_n)$ be a decreasing
sequence of closed bounded subsets of $\ell^2(X)$ whose diameters go to zero.  For each $x\in X$,
let $A_n^x$ be the set of $x$'th coordinates of elements of $A_n$, i.e.\ $A_n^x=P_x(A_n)$, where
$P_x$ is the orthogonal projection of $ \calH$ onto the one-dimensional subspace spanned by $\xi_x$.
Then $(\overline{A_n^x})$ is a decreasing sequence of closed bounded subsets of $\bbC$ whose diameters
go to 0; hence $\cap_n \overline{A_n^x}=\{c_x\}$ for some $c_x\in\bbC$.  We claim
$$\sum_{x\in X}c_x\xi_x\in\bigcap_n A_n\ .$$
If $F$ is a finite subset of $X$, set $A_n^F=P_F(A_n)$, where $P_F$ is the orthogonal projection
of $\ell^2(X)$ onto the span $\cY_F$ of $\{\xi_x:x\in F\}$.  Then $(\overline{A_n^F})$ is a decreasing sequence
of closed subsets of $\cY_F$ whose diameters go to zero, so $\cap_n \overline{A_n^F}=\{\eta_F\}$ for some
$\eta_F\in\cY_F$ since $\cY_F$ is finite-dimensional; we must have $\eta_F=\sum_{x\in F}c_x\xi_x$.

\smallskip

\noindent
Suppose $\|\zeta\|\leq M$ for all $\zeta\in A_1$.  It follows that for each finite subset $F$ of $X$,
$\|\zeta\|\leq M$ for all $\zeta\in A_1^F$; so $\sum_{x\in F}|c_x|^2\leq M^2$ for any $F$.
Thus $\sum_{x\in X}|c_x|^2\leq M^2$, so the $c_x$ define a vector
$$\eta=\sum_{x\in X} c_x\xi_x\in\ell^2(X)\ .$$

\smallskip

\noindent
%To show $\eta\in\cap_n A_n$, define a directed set
%$$D=\left \{ (\zeta,F,n): F\subseteq X\mbox{ finite, }n\in\N,\zeta\in A_n, \|P_F(\zeta)-\eta_F\|<\frac{1}{n}\right \}$$
%with preordering $(\zeta,F,n)\leq(\theta,E,m)$ if $F\subseteq E$ and $n\leq m$.  Define a net
%in $ \calH$ by $\theta_{(\zeta,F,n)}=\zeta$.  Then $\theta_{(\zeta,F,n)}\to\eta$.  Since
%$\theta_{(\zeta,F,n)}\in A_n$ and each $A_n$ is closed, $\eta\in\cap_n A_n$.
%
To show $\eta\in\cap_nA_n$, let $\epsilon>0$, and fix $m$ with $\diam(A_m)\leq\epsilon$.  Let $n\geq m$.
If $F$ is a finite subset of $X$, we have $\diam(\overline{A_n^F})=\diam(A_n^F)\leq\epsilon$, so, for any $\zeta\in A_n$,
$$\|P_F(\zeta)-\eta_F\|^2=\sum_{x\in F}|\langle\zeta,\xi_x\rangle-\langle\eta,\xi_x\rangle|^2
=\sum_{x\in F}|\langle\zeta-\eta,\xi_x\rangle|^2\leq\epsilon^2\ .$$
Thus $\|\zeta-\eta\|\leq\epsilon$ for all $\zeta\in A_n$.  This is true for every $\epsilon>0$, for all
sufficiently large $n$.  Since the $A_n$ are decreasing and closed, $\eta\in\cap_nA_n$.
\end{Proof}

\paragraph{}
If $X$ is a set, $\ell^2(X)$ has a canonical orthonormal basis $\{\xi_x:x\in X\}$, where $\xi_x$
is the characteristic (indicator) function of $\{x\}$.  If $\eta\in\ell^2(X)$, then we have
$$\eta=\sum_{x\in X} \eta(x)\xi_x$$
in the sense that the sums over finite subsets converge in norm to $\eta$.  Every square-summable series occurs: if $\{c_x:x\in X\}$ is a square-summable set of complex
numbers, there is a unique $\eta\in\ell^2(X)$ with $c_x=\langle\eta,\xi_x\rangle$ for all $x\in X$.

\bigskip

Conversely:

\paragraph{}
\begin{Proposition}\label{OrthBasisL2}
Let $ \calH$ be a Hilbert space with an orthonormal basis $\{\zeta_x:x\in X\}$ for some set $X$.
Then the map $\zeta_x\mapsto\xi_x$ gives an isometric isomorphism from $ \calH$ onto $\ell^2(X)$.
\end{Proposition}

\begin{Proof}
Let $ \calH_0$ be the dense subspace of $ \calH$ of finite linear combinations of the $\zeta_x$.
The map $\zeta_x\mapsto\xi_x$ defines an isometric linear map of $ \calH_0$ onto a dense subspace
of $\ell^2(X)$.  Since
$\ell^2(X)$ is $\sigma$-complete, this map extends (uniquely) to an isometric linear map $T$ from
$ \calH$ onto $\ell^2(X)$ by \ref{DenseIsom}.
\end{Proof}

\smallskip

Thus the $\ell^2(X)$ for various $X$ are universal models of Hilbert spaces with orthonormal bases.

\section{Finite-Dimensional and Separable Hilbert Spaces}\label{S.FiniteSeparable}

\paragraph{}\label{FinDimHilbSp}
No Choice is needed to prove the standard results of finite-dimensional linear algebra.  If $ \calH$
is a finite-dimensional inner product space (finite-dimensional means spanned by a finite set
of vectors), then:
\begin{enumerate}
\item[(i)]  Every set of orthonormal vectors can be expanded to an orthonormal basis. In
particular, $ \calH$ has an orthonormal basis.
\item[(ii)]  $ \calH$ has a well-defined orthogonal dimension $n$, which coincides with its linear dimension.
\item[(iii)]  Every linearly independent set can be orthonormalized (Gram-Schmidt).
\end{enumerate}
In addition, we have (in ZF):
\begin{enumerate}
\item[(iv)]  $ \calH$ is $\sigma$-complete.
\item[(v)]  The closed unit ball of $ \calH$ is compact.
\end{enumerate}

\paragraph{}\label{InfDimTotBdd}
There is a converse to (v):  if $ \calH$ is an infinite-dimensional inner product space,
then for any $n$ there is an orthonormal set of vectors in $ \calH$ of cardinality $n+1$ made by Gram-Schmidt, so the closed
unit ball cannot be covered by $n$ open balls of radius $\frac{1}{\sqrt{2}}$ or less, and thus the closed
unit ball of $ \calH$ is not totally bounded.

(No Choice is needed to show that a compact metric space is $\sigma$-complete and totally bounded,
but the converse requires some Choice, cf.\ \cite[XII.16.12.9]{BlackadarReal}.)

\subsection{Separable Hilbert Spaces}

\paragraph{}
Most of the results about finite-dimensional Hilbert spaces extend to separable Hilbert spaces
(every finite-dimensional Hilbert space is separable),
due to the Gram-Schmidt orthogonalization process where any sequence of vectors in an inner
product space can be orthogonalized into an orthonormal sequence (or finite sequence) with the
same span.  The procedure is inductive (recursive) and requires no Choice and, indeed, no completeness.
%See e.g.\ \cite[]{BlackadarReal} for details of the process.
Using Gram-Schmidt plus previous results, we have:

\paragraph{}
\begin{Theorem}\label{SepHilbSpaceThm}
Let $ \calH$ be a separable Hilbert space.  Then
\begin{enumerate}
\item[(i)]  $ \calH$ has a countable orthonormal basis.
\item[(ii)]  Every closed subspace of $ \calH$ is separable (hence has a countable orthonormal basis).
\item[(iii)]  Every orthonormal set in $ \calH$ can be expanded into an orthonormal basis.
\item[(iv)]  Any two orthonormal bases for $ \calH$ have the same cardinality (finite or $\aleph_0$).
\end{enumerate}
Conversely, any Hilbert space with a countable orthonormal basis is separable.
\end{Theorem}

\begin{Proof}
(i) is exactly the Gram-Schmidt process.  For (ii), if $\cY$ is a closed subspace of $ \calH$, then there
is a projection of norm one from $ \calH$ onto $\cY$ (\ref{OrthDecompProp}), and a continuous image
of a separable space is separable.  For (iii), let $\cS=\{\xi_x:x\in X\}$ be an orthonormal set in $ \calH$,
and let $\cY$ be the closed linear span of $\cS$.  Then $\cY^\perp$ is a separable Hilbert space,
hence has a (countable) orthonormal basis $\cB$, and $\cB\cup\cS$ is an orthonormal basis
for $ \calH$.  For the converse, finite linear combinations of the basis vectors with coefficients in
$\Q+\Q i$ are dense, and can be effectively enumerated without any Choice.

\smallskip

\noindent
(iv)  Let $ \calH$ be a separable Hilbert space, and $\{\xi_x:x\in X\}$ be an orthonormal set in $ \calH$.
Let $\{\eta_n:n\in\N\}$ be a countable dense set in $ \calH$.  For each $x\in X$, let $n(x)$ be the
smallest $n$ for which $\|\eta_n-\xi_x\|<\frac{1}{2}$.  Then $x\to n(x)$ is injective, so $X$ is
countable.  If $ \calH$ is finite-dimensional, an orthonormal basis for $ \calH$ is a Hamel (vector space)
basis for $ \calH$, and the vector space dimension of a finite-dimensional vector space is well
defined by elementary linear algebra (no Choice needed).  If $ \calH$ is infinite-dimensional, then
every orthonormal basis for $ \calH$ is countably infinite.
\end{Proof}

Theorem~\ref{SepHilbSpaceThm} is about the existence of a well-ordered basis rather than separability.

\paragraph{}\begin{Theorem} \label{L.l2(kappa)} Suppose that $\kappa$ is a well-ordered cardinal.
	Let $ \calH$ be a Hilbert space with a dense subset of cardinality $\kappa$.  Then
	\begin{enumerate}
		\item[(i)]  $ \calH$ has an orthonormal basis of cardinality $\lambda$, for some $\lambda\leq \kappa$.
		\item[(ii)]  Every closed subspace of $ \calH$ has orthonormal basis of cardinality $\leq \kappa$.
		\item[(iii)]  Every orthonormal set in $ \calH$ can be expanded into an orthonormal basis.
		\item[(iv)]  Any two orthonormal bases for $ \calH$ have the same cardinality.
	\end{enumerate}
\end{Theorem}

The proof of this theorem is given below after a definition and a lemma.

\paragraph{}\begin{Definition}The following describes the \emph{transfinite Gram--Schmidt process}.
	Suppose that $\kappa$ is an ordinal and $\eta_\alpha$, for $\alpha<\kappa$, are nonzero vectors in a Hilbert space. Define vectors $\zeta_\alpha$, for $\alpha<\kappa$,  by transfinite recursion as follows.

	First let $\zeta_0=\eta_0\|\eta_0\|^{-1}$. If $\zeta_\alpha$, for $\alpha<\beta$, have been determined, then consider the projection of $\xi_\beta$ to the closed linear span of $\zeta_\alpha$, for $\alpha<\beta$:
	\begin{equation}\label{eq.GramSchmidt}
		\xi_\beta=\sum_{\alpha<\beta} \langle \eta_\beta, \zeta_\alpha\rangle \zeta_\alpha.
	\end{equation}
	If $\eta_\beta-\xi_\beta$ is nonzero, let $\zeta_\beta=(\eta_\beta-\xi_\beta)\|\eta_\beta-\zeta_\beta\|^{-1}$. Otherwise let $\zeta_\beta=0$.
\end{Definition}

\paragraph{}\begin{Lemma}\label{L.transfinite.GS}
	If $\zeta_\alpha$, for $\alpha<\kappa$, are vectors in a Hilbert space indexed by an ordinal, then the transfinite Gram--Schmidt process results in a well-ordered orthonormal basis for the closed linear span of these vectors.
\end{Lemma}

\begin{Proof} The proof is analogous to one in the finitary case. At the limit stages, Bessel’s inequality implies that the finite partial sums of the right-hand side of \eqref{eq.GramSchmidt} converge.  In the resulting sequence $(\zeta_\alpha)$, the non-zero vectors are well-ordered by some ordinal $\kappa'\leq \kappa$.
\end{Proof}

Notably, the use of completeness at limit stages of the transfinite Gram--Schmidt process is necessary even if the   AC holds, if $H$ is a non-separable (not necessarily $\sigma$-complete) inner product space then it may not contain an orthonormal basis (see \cite{DixmierBases}).

\begin{ProofOf}{Proof of Theorem~\ref{L.l2(kappa)}} The proof is analogous to the proof of Theorem~\ref{SepHilbSpaceThm}.

	(i) Let $\eta_\alpha$, for $\alpha<\kappa$, be an enumeration of a dense subset of $ \calH$. By applying the transfinite Gram--Schmidt process and using Lemma~\ref{L.transfinite.GS} we obtain an ordinal $\lambda\leq \kappa$ such that $ \calH\cong \ell^2(\lambda)$. By re-enumerating, we can assume $\lambda$ is a cardinal.

	(ii) If $ \calH_0$ is a closed subspace of $ \calH$, consider the image of a dense subset of $ \calH$ under the projection to $ \calH_0$ and apply (i).

	(iii) Apply (ii) to the orthogonal complement of the span of the given orthonormal set.

	(iv) By (ii), we need to prove that if $\kappa$ and $\lambda$ are well-ordered cardinals and $\ell^2(\kappa)$ and $\ell^2(\lambda)$ are isomorphic, then $\kappa=\lambda$.  This is elementary
if $\kappa$ is finite.  For the infinite case, by expanding each
basis vector in terms of the other basis and using that the support of a vector is countable (\ref{VectSuppC}, and we are working in a well-ordered set), we
obtain $\kappa\leq\aleph_0\cdot\lambda$ and $\lambda\leq\aleph_0\cdot\kappa$.  But if $\kappa$
is a well-ordered infinite cardinal, then $\kappa=\aleph_0\cdot\kappa$.
\end{ProofOf}

\section{Dedekind-Finite Sets of Various Flavors}\setcounter{paragraph}{0}
\label{S.DF+}
In the absence of CC, there can be bizarre sets which are infinite but ``almost finite.''
In this section we survey various types of such sets which exist in some models of ZF.
For more details see e.g.\ \cite{jech2008axiom}, also \S\ref{S.ModelsofZF}.

We start by making precise the notion of cardinality in the context of ZF.

\subsection{Vanilla Flavor}\label{S.DF}
\paragraph{}\begin{Definition}\label{Def.Cardinal}
	The \emph{cardinality} of a set $X$, denoted $|X|$, is the equivalence class of all $Y$ such that there is a bijection between $X$ and $Y$. (This is a proper class, but by Scott’s trick it can be identified with a set.) AC is equivalent to the assertion that every set can be well-ordered, and one writes $\kappa=|X|$ where $\kappa$ is the least ordinal for which such bijection exists. We take such $\kappa$ (if it exists) as the canonical representative of $|X|$, and refer to it as a \emph{well-ordered cardinal}.
\end{Definition}

We will use symbols $\kappa,\lambda$ to denote cardinals, with understanding that they are not necessarily well-orderable.

The cardinality of $X$ is less or equal than the cardinality of $Y$, in symbols $|X|\leq |Y|$, if there is an injection from $X$ to $Y$.
If $Y$ is nonempty, then $|Y|\leq |X|$ implies that there is a surjection from $X$ onto $Y$. (To see this, fix an injection $f\colon Y\to X$ and $y_0\in Y$. Define $g\colon X\to Y$ by $g(f(y))=y$ and $g(x)=y_0$ if $x$ is not in the range of $g$.)

\paragraph{}
We use the following notation for a set $X$.  Denote by  $\cP(X)$
the power set of $X$, the set of all subsets of $X$.  Let $ \Fin(X)$ be the subset of $\cP(X)$
consisting of all finite subsets of $X$ (set theorists often denote this by $[X]^{<\omega}$).
$\Seq(X)$ is the set of all finite sequences of distinct elements of $X$ (including the ``empty sequence'').

We have $|X|\leq| \Fin(X)|\leq|\cP(X)|$ and, of course, $|X|<|\cP(X)|=2^{|X|}$; and $|X|\leq|\Seq(X)|$.
There is an obvious surjective map from $\Seq(X)$ to $ \Fin(X)$.  We will have use for the
following well-known fact, a theorem of ZF (i.e.\ no Choice is needed):

\paragraph{}
\begin{Proposition}\label{CountFinProp}
Let $\N$ be the set of natural numbers.  Then $\Seq(\N)$ and $ \Fin(\N)$ are countably infinite
(equipotent with $\N$).
\end{Proposition}
\paragraph{}
We write $|Y|\leq^*|X|$ if there is a surjective function from $X$ to $Y$ (or if $Y=\emptyset$).
We have $|Y|\leq|X| \imp |Y|\leq	^*|X|$, and the converse is true under AC but not otherwise (whether this is equivalent to AC is an old open problem).

\bigskip

The following well-known notion is central to our study.
\paragraph{}
\begin{Definition}\label{DedekindFinSetDef}
A set  $X$ is {\em Dedekind-finite} if there is no bijection from $X$ to a proper subset of~$X$, i.e.\ every
injective function from $X$ to $X$ is surjective.  Otherwise $X$ is {\em Dedekind-infinite}. A set is \emph{finite} if its cardinality is equal to the cardinality of some natural number.
A cardinal $\kappa$ is Dedekind-finite if it is the
cardinal of a Dedekind-finite set.  Say a set $X$ or a cardinal $\kappa$ is {\em DF} if it is infinite and Dedekind-finite.
\end{Definition}

\smallskip

 While every finite set is Dedekind-finite, the converse fails in some models of ZF (Proposition~\ref{P.Cohen}).
Whether a set $X$ is Dedekind-finite or Dedekind-infinite depends only on the cardinal of $X$:
if $X$ and $Y$ are equipotent and one is Dedekind-finite, so is the other.  Any subset of a
Dedekind-finite set is Dedekind-finite.

The next result is fundamental:

\paragraph{}
\begin{Proposition}\label{DedFinProp}
Let $X$ be a set.  Then $X$ is Dedekind-infinite if and only if there is an injective function from $\N$
to $X$ (i.e.\ there is a sequence of distinct elements of $X$).  A cardinal $\kappa$ is
Dedekind-infinite if and only if $\aleph_0\leq\kappa$.
\end{Proposition}

\smallskip

\begin{Proof} If $|Y|\leq |X|$ and  $Y$ is Dedekind-infinite, then so is $X$; in particular, if there is an injection from $\bbN$ into $X$ then $X$ is Dedekind-infinite. Suppose that $X$ is Dedekind-infinite and let $f\colon X\to X$ be an injection that is not a surjection. Choose $x\in X\setminus f[X]$.  Writing $f^n$ for the $n$-fold iteration of $f$ and $f^0$ for the identity map,  $g\colon \bbN\to X$ defined by $g(n)=f^n(x)$ is an injection from $\bbN$ into $X$.

	The second sentence is a restatement of the first.
\end{Proof}

\paragraph{}
\begin{Corollary} \label{C.DF+finite} If $X$ is a proper subset of $Y$ and $X$ is DF, then $|X|<|Y|$.
\end{Corollary}

If the CC is assumed, a set is Dedekind-infinite if and only if it is infinite, so a set is
Dedekind-finite if and only if it is finite; thus DF sets are incompatible with the CC.

\bigskip

The next proposition gives some alternate characterizations of Dedekind-infiniteness:

\paragraph{}
\begin{Proposition}
Let $X$ be a set, and $\kappa=|X|$.  The following are equivalent:
\begin{enumerate}
\item[(i)]  $X$ is Dedekind-infinite, i.e.\ there is a proper subset $Y$ of $X$ which is equipotent with $X$.
\item[(ii)]  $\kappa=\kappa+\aleph_0$.
\item[(iii)]  $\kappa=\kappa+1$.
\item[(iv)]  If $\kappa=\lambda+1$, then $\lambda=\kappa$.
\item[(v)]  If $Y$ is any set obtained from $X$ by adding or removing finitely many points, then $|Y|=\kappa$.
%\item[(vi)]  There is a proper subset $Y$ of $X$ which is equipotent with $X$.
\end{enumerate}
\end{Proposition}

%See () for the proof.

\paragraph{}\label{DedFinSubsets}
If we have one DF cardinal, we have many.  If $X$ is a DF
set and $n\in\N$, let $F_n$ be an $n$-element subset of $X$, and $\kappa_n=|X\setminus F_n|$
(note that we can choose $F_n$ for any particular $n$, or for finitely many $n$, but we may not
choose $F_n$ simultaneously for all $n$, and in particular we cannot choose the $F_n$ so that
$F_n\subseteq F_{n+1}$ for all $n$).  It is easily seen that the cardinality of $X\setminus F_n$
depends only on $n$ and not on the choice of $F_n$, so $\kappa_n$ is well defined for each $n$.
We clearly have $\kappa_{n+1}\leq\kappa_n$ for each $n$, and we have $\kappa_{n+1}\neq\kappa_n$
since $X\setminus F_n$ is DF.  Thus there is a strictly decreasing sequence of DF
cardinals smaller than $|X|$, which dramatically contradicts the Well-Ordering
Principle (there is even a collection of DF cardinals order-isomorphic to $\R$,
cf.\ \cite[II.9.5.21]{BlackadarReal}).  More generally, if $Y$ is any proper subset of $X$, then $|Y|<|X|$.  DF cardinals do not
even form a set in general: in fact, there are models of ZF in which every set is the image of a DF set. This holds in the `improved Morris model'  of \cite[\S 5.8]{karagila2019iterated}, and we briefly explain why.  It suffices to prove that every rank-initial segment of the universe, $V_\alpha$ is the image of a Dedekind-finite set.
  \cite[Proposition 4.2]{karagila2019iterated} shows that the set of generic branches is Dedekind-finite, and by construction it can be mapped onto $V_\alpha$; as the forcing described in \cite{monro1975independence} is iterated through the universe, in the final model each $V_\alpha$ is the image of some Dedekind-finite set.

\paragraph{}
\begin{Proposition}\label{DFProdProp}
A finite Cartesian product of Dedekind-finite sets is Dedekind-finite.  A finite union of
Dedekind-finite sets is Dedekind-finite.
\end{Proposition}

\begin{Proof}
If $(x_j)$ is a sequence of distinct elements of $X_1\times\cdots\times X_n$, then for some $k$,
$1\leq k\leq n$, the sequence of $k$'th coordinates of the $x_j$ must contain infinitely many
distinct elements of $X_k$, hence $X_k$ is not Dedekind-finite.  Similarly, if $X=\cup_{k=1}^nX_k$,
and $(x_j)$ is a sequence of distinct elements of $X$, then infinitely many of the $x_j$ must be
in some $X_k$, so $X_k$ is not Dedekind-finite.
\end{Proof}

\smallskip
In particular, if $X$ is Dedekind-finite, then $X^n$ and $X\times F$ are Dedekind-finite for any
$n$ and any finite set $F$.  These sets all have different cardinalities by \ref{DedFinSubsets}.
By repeating these constructions and taking subsets, a very large collection of Dedekind-finite
cardinals can be obtained.

\subsection{Other Flavors}

Here are some important variations on Dedekind-infiniteness:

\paragraph{}
\begin{Definition}\label{WeakDedInfDef}
Let $X$ be a set. Then
\begin{enumerate}
\item[(i)]  $X$ is {\em Cohen-infinite} if $ \Fin(X)$ is Dedekind-infinite, i.e.\ there is a sequence
of distinct finite subsets of $X$.  %$X$ is {\em Cohen-finite} if it is not Cohen-infinite.
\item[(ii)]  $X$ is {\em weakly Dedekind-infinite} $\cP(X)$ is Dedekind-infinite, i.e.\ there is a sequence
of distinct subsets of $X$.  %$X$ is {\em power Dedekind-finite} if it is not weakly Dedekind-infinite.
\item[(iii)]  $X$ is {\em dually Dedekind-infinite} if there is a surjective function from $X$ to $X$
which is not injective.  %$X$ is {\em dually Dedekind-finite} if it is not dually Dedekind-infinite.
\end{enumerate}
\end{Definition}

\paragraph{}
We have ($X$ Dedekind-infinite) $\imp$ ($X$ Cohen-infinite) $\imp$ (X weakly Dedekind-infinite) $\imp$ (X infinite) and ($X$ Dedekind-infinite) $\imp$ ($X$ dually Dedekind-infinite).

\paragraph{}
\begin{Definition}\label{Def.CF}
Let $X$ be a set.  Then
\begin{enumerate}
%\item[(i)]  $X$ is {\em dually Dedekind-finite} if it is not dually Dedekind-infinite.
\item[(i)]  $X$ is {\em Cohen-finite} if $ \Fin(X)$ is
Dedekind-finite.
\item[(ii)]  $X$ is {\em power Dedekind-finite} if $\cP(X)$ is
Dedekind-finite, i.e.\ if $X$ is not weakly Dedekind-infinite.
\item[(iii)]  $X$ is {\em dually Dedekind-finite} if it is not dually Dedekind-infinite.
\end{enumerate}
\end{Definition}

\smallskip

We use the following abbreviations:
\begin{enumerate}
\item[]  CF means ``infinite and Cohen-finite.''
\item[]  PF means ``infinite and power Dedekind-finite.''
\item[]  DDF means ``infinite and dually Dedekind-finite.''
\end{enumerate}
\smallskip

We use the term {\em power Dedekind-finite} instead of {\em weakly Dedekind-finite}
for the negation of {\em weakly Dedekind-infinite}
since the condition of being power Dedekind-finite is considerably {\em stronger} than being
Dedekind-finite (\ref{DFnotPF}).

The terminology {\em Cohen-finite} is justified since by \cite[\S 6]{karagila2020have}, a set $X$ is CF if and only if the forcing $\Add(X,1)$ for adding a subset of $X$ by finite conditions (`adding a Cohen subset of $X$') has the property that every  statement in the language of forcing is decided by a finite predense set. This property is also equivalent to the assertion that  $X$ remains DF in the forcing extension by $\Add(X,1)$.

Any subset of a Cohen-finite [resp.\ power Dedekind-finite, dually Dedekind-finite] is
Cohen-finite [resp.\ power Dedekind-finite, dually Dedekind-finite].
The definitions easily imply that ($X$ PF) $\imp$ ($X$ CF) $\imp$ (X DF) and ($X$ DDF) $\imp$ ($X$ DF), and  \ref{PFFns} implies ($X$ PF) $\imp$ ($X$ DDF).  None of these implications can be reversed (\ref{DFnotDDF} shows that DF sets are not necessarily DDF, \ref{RussellNotCF} shows that some DF sets are not CF, and combining \ref{DedInfRSubProp} and \ref{CohenModel} shows that CF sets are not always PF).
In addition to that, ($X$ DDF) does not imply ($X$ CF), as shown in \ref{RussellDDF}, nor ($X$ CF) implies ($X$ DDF), as shown by \ref{DFnotDDF} when considered in Cohen's model.

\bigskip

The next result, due to {\sc Kuratowski} \cite[pp.94--95]{TarskiFinis}, gives alternate characterizations of weak Dedekind-infiniteness.  See e.g.\ \cite[II.9.5.7]{BlackadarReal} for a proof.

\paragraph{}
\begin{Theorem}\label{WeakDedInfThm}
Let $\kappa$ be a cardinal, and $X$ a set with $|X|=\kappa$.  The following are equivalent:
\begin{enumerate}
\item[(i)]  $\aleph_0\leq^*\kappa$, i.e.\ there is a sequence of pairwise disjoint nonempty subsets of $X$.
\item[(ii)]  $\aleph_0\leq 2^\kappa$, i.e.\ $X$ is weakly Dedekind-infinite.
\item[(iii)]  $2^{\aleph_0}\leq 2^\kappa$.
\end{enumerate}
\end{Theorem}

\paragraph{}
\begin{Proposition}\label{DedInfRSubProp}
	Let $X$ be an infinite second-countable Hausdorff space (e.g.\ an infinite subset of $\R$).  Then $X$ is weakly Dedekind-infinite.
\end{Proposition}

\begin{Proof}
	We will prove that $X$ maps onto $\bbN$.  Fix an enumeration $U_n$, for $n\in \bbN$, of a base for the topology of $X$ with $U_0=X$. Suppose for a moment that all points of $X$ are isolated. Then $f(x)=\min\{n\mid U_n\cap X=\{x\}\}$ defines a surjection from $X$ onto an infinite subset of $\bbN$.

	We may therefore assume that $X$ has an  accumulation point, $x$. Recursively define $m(n)\in \bbN$, for $n$ in $\bbN$, by $m(0)=0$, and if $m(n)$ has been defined let
	\[
	m(n+1)=\min\{j \mid U_j\subseteq U_{m(n)},\text{ $x\in U_j$,  and } (U_{m(n)}\setminus U_j)\cap X\neq \emptyset\}.
	\]
	Then the function that sends $X\cap (U_{m(n)}\setminus U_{m(n+1)})$ to $n$ is a surjection.
\end{Proof}
\smallskip

See \cite{SierpinskiCardinal} or \cite[II.9.5.6]{BlackadarReal} for another proof.  Subsets of $\R$
are special among all sets (in the absence of AC) in part because they can be totally ordered.

\paragraph{}
\begin{Proposition}\label{StrongDFSubsets}
	Let $X$ be a CF set.  Then every subset of $X$ that can be covered by a countable union of finite sets is finite.
\end{Proposition}

\begin{Proof}
	Let $A=\cup_{n=1}^\infty A_n\subseteq X$, with each $A_n$ finite.  If $A$ is infinite, then there must be
	infinitely many distinct $A_n$, hence a sequence of distinct $A_n$, so $ \Fin(X)$ is not Dedekind-finite.
\end{Proof}

\paragraph{}
The analog of \ref{WeakDedInfThm} (ii) $\imp$ (i) (which is quite subtle) for finite subsets is easily
proved by induction.  Thus a set is Cohen-finite if and only it it does not contain a sequence of
pairwise disjoint nonempty finite subsets.

\paragraph{}
\begin{Proposition}\label{P.lo.CF}
	If $X$ can be linearly ordered, then it is DF if and only if it is CF.
\end{Proposition}

\begin{Proof}
	Only the direct implication requires a proof, and we will prove the contrapositive. Assume $X$ is linearly ordered by $\prec$ and $F(n)$ is a sequence of distinct nonempty finite subsets of $X$. By the pigeonhole principle, $\bigcup_n F(n)$ is infinite. We can order $\bigcup_n F(n)$ by (writing $n(x)=\min\{n: x\in F(n)\}$)
	\[
	x\prec_1 y\text{ if and only if $n(x)<n(y)$ or $n(x)=n(y)$ and $x\prec y$}.
	\]
	Since each $F(n)$ is finite, this is a well-ordering of the infinite subset $\bigcup_n F(n)$ of $X$. Therefore $X$ is not Dedekind-finite.
\end{Proof}

\bigskip

We have the following analog of \ref{DFProdProp}:

\paragraph{}
\begin{Proposition}\label{PFProdProp}
(i)  A finite Cartesian product of Cohen-finite sets is Cohen-finite.  A finite union of
Cohen-finite sets is Cohen-finite.

\smallskip

\noindent
(ii)  A finite Cartesian product of power Dedekind-finite sets is power Dedekind-finite.  A finite union of
power Dedekind-finite sets is power Dedekind-finite.
\end{Proposition}

\begin{Proof}
The proof of (i), and (ii) for unions, is straightforward along the lines of the proof of \ref{DFProdProp}.
For products of power Dedekind-finite sets, it suffices to show the result for the Cartesian product
of two power Dedekind-finite sets $X$ and $Y$.  Suppose $f$ is a surjective function from
$X\times Y$ to $\N$.  For each $x\in X$, the restriction of $f$ to $\{x\}\times Y$ has finite range $A_x$
since $\{x\}\times Y$ is power Dedekind-finite.  Thus $x\mapsto A_x$ is a function from $X$ to
$ \Fin(\N)$, which must have infinite range since the union of the $A_x$ is $\N$.  But $ \Fin(\N)$
is countable (\ref{CountFinProp}), contradicting that $X$ is power Dedekind-finite.
\end{Proof}

\paragraph{}\label{PFFns}
An interesting consequence is that if $X$ and $Y$ are PF, then there is not a sequence of
distinct relations from $X$ to $Y$, and in particular there is no sequence of distinct functions from
$X$ to $Y$.

If $X$ is dually Dedekind-infinite, and $f:X\to X$ is surjective but not injective, then $(f^n)$
is a sequence of distinct functions from $X$ to $X$, and hence $X$ cannot be power Dedekind-finite.
Thus $X$ PF $\imp$ $X$ DDF.

\paragraph{}\label{DFnotDDF}
If $X$ is DF, then $\Seq(X)$ is also DF [a sequence of distinct elements of $\Seq(X)$ can be concatenated and the pigeonhole principle implies that it contains a subsequence of distinct elements of $X$].  However, if $X$ is
infinite, $\Seq(X)$ is not DDF, since there is a surjective and noninjective function from
$\Seq(X)$ to itself sending the empty sequence to itself and dropping the first term in each nonempty
sequence.  Thus, if there exists a DF set, there exists a DF set which is not DDF.

\paragraph{}
The situation with $ \Fin(X)$ is quite different, even though there is a finite-to-one surjective function
from $\Seq(X)$ to $ \Fin(X)$.  If $X$ is DF, then $ \Fin(X)$ can be DF, even DDF, or not DF.
See \cite{karagila2020have} for a study.  We do have the following (for $n\geq 1$ let $\Fin^n(X)$ be  the $n$-th iterate of the operation $\Fin$).

\paragraph{} \begin{Proposition}\label{P.CFCF} If $X$ is a CF set, then  $\Fin^n(X)$ is CF for all $n\geq 1$.
\end{Proposition}

\begin{Proof} Suppose $X$ is CF. In order to prove that $ \Fin(X)$ is CF,  towards contradiction let  $F_n$ be an infinite sequence of disjoint nonempty subsets of $ \Fin(X)$.
	Let $G_n=\bigcup F_n$. Then $G_n$ are finite subsets of $X$, and since $X$ is CF there is $G\in  \Fin(X)$ such that $G_n=G$ for infinitely many $n$. However, both $\cP(G)$ and $\cP(\cP(G))$ are finite, contradicting the assumption that $F_n$ were pairwise distinct.

	This proves that $ \Fin(X)$ is CF, and it also provides the inductive step for the proof that $\Fin^n(X)$ is CF for all $n$.
\end{Proof}

\paragraph{}\label{DFnotPF}
If $X$ is infinite, then $\cP( \Fin(X))$ (and {\em a fortiori} $\cP(\cP(X))$) is always Dedekind-infinite:
the sets $A_n$ of $n$-element subsets of $X$ form a sequence of distinct elements in $\cP( \Fin(X))$.
However, if $X$ is a DF subset of $\R$, then $X$ is CF (Proposition~\ref{P.lo.CF}) and therefore $ \Fin(X)$ is DF and even CF (\ref{P.CFCF}), although $X$ is not PF (\ref{DedInfRSubProp}).

\subsection{Amorphous Sets}

We now describe some particularly strange DF sets whose existence cannot be refuted in ZF set theory.\footnote{More precisely, if ZF has a model, then it has a model with an amorphous set. This also applies to other flavors of DF sets discussed in \S\ref{S.DF} and \S\ref{S.DF+}. The reason why it is necessary to assume that ZF has a model is that  By G\"odel's Incompleteness Theorem, in ZF there is no proof that there is a model of ZF\dots{} unless ZF is inconsistent.}

\paragraph{}
\begin{Definition}\label{AmorphSetDef}
An infinite set $X$ is {\em amorphous} if $X$ cannot be written as a disjoint union of two
infinite subsets, i.e.\ every subset of $X$ is either finite or cofinite.
\end{Definition}

\paragraph{}\label{AmorphSetDedFin}
Any infinite subset of an amorphous set is amorphous.
If $X$ is an amorphous set, $Y$ an infinite set, and there is a surjective function $f:X\to Y$, then $Y$ is
also amorphous [if $Y$ is the disjoint union of $Y_1$ and $Y_2$, then $X$ is the disjoint union
of $f^{-1}(Y_1)$ and $f^{-1}(Y_2)$].  In particular, if $X$ is weakly Dedekind-infinite, since $\N$
is not amorphous~$X$ cannot be amorphous, so if~$X$ is amorphous then $X$ is PF.
Thus amorphous sets are Dedekind-finite in a strong sense.  In particular, amorphous
sets are incompatible with the CC.

\smallskip

As an example of how strange an amorphous set is, we have:

\paragraph{}
\begin{Proposition}\label{AmorphSetProp}
An amorphous set cannot be totally ordered.
\end{Proposition}

\smallskip

See e.g.\ \cite [II.9.5.13]{BlackadarReal} for a proof.
\subsection{Permanence Properties of CF Sets}
\label{S.Permanence}

\paragraph{} For a relation $R\subseteq X\times Y$ consider its horizontal and vertical sections,  $R_x=\{y\in Y: (x,y)\in R\}$ and $R^y=\{x\in X: (x,y)\in R\}$.

\paragraph{}\begin{Lemma} \label{L.BackAndForth} Suppose that $X$ and $Y$ are sets and at least one of them  is CF.  Also suppose that $R\subseteq X\times Y$, and all horizontal and vertical sections of $R$ are finite.
	Then there are $Z$ and partitions $X=\bigsqcup_{z\in Z} X_z$, $Y=\bigsqcup_{z\in Z} Y_z$ into finite sets such that $\bigcup\{R_x: x\in X_z\}\subseteq Y_z$  and $\bigcup\{R^y: y\in Y_z\}\subseteq X_z$ for all $z\in Z$.
	This implies that both sets are  CF.
\end{Lemma}

\begin{Proof}  Note that all sections of $R$ are finite if and only if all sections of the inverse relation $R^{-1}$ are finite. It therefore suffices to prove the case when  $X$ is CF, since otherwise we can exchange the roles of $X$ and $Y$ and consider $R^{-1}$ instead of $R$.

	Fix for a moment $x\in X$.
	Recursively define $F(n)\subseteq X$ and $G(n)\subseteq Y$ as follows. Let $F(0)=\{x\}$, and for all $n$ let
	\begin{align*}
		G(n)&=\bigcup\{R_{x'}: x'\in F(n)\},\\
		F(n+1)&= \bigcup\{R^y: y\in G(n)\}.
	\end{align*}
	Clearly $F(n)\subseteq F(n+1)$ and $G(n)\subseteq G(n+1)$ for all $n$.  By induction on $n$ one proves that all $F(n)$ and all $G(n)$ are finite, being unions of finitely many finite sets. Since $X$ is CF, there exists $m=m(x)$ such that $F(m)=F(n)$ for all $n\geq m$. This implies that $G(m)=G(n)$ for all $n\geq m$.
	Let $F[x]=F(m)$ and $G[x]=G(m)$.

	Fix $x'\in F[x]$. By induction on the minimal $n$ such that $x'\in F(n)$, one proves that $F[x']=F[x]$, and therefore $G[x']=G[x]$. Let $Z=\{F[x]: x\in X\}$. For every $z\in Z$ there is a unique finite subset $G[z]$ of $Y$ which satisfies $z=G[x]$ for some (every) $x$ such that $F[x]=z$.  Then $X_z=z$ and $Y_z=G[z]$ define  partitions of $X$ and $Y$ into finite sets.  By construction, $\bigcup\{R_x: x\in X_z\}\subseteq Y_z$  and $\bigcup\{R^y: y\in Y_z\}\subseteq X_z$ for all $z\in Z$.

	To see that $Y$ is CF, fix an increasing sequence of finite subsets of $Y$, $G(n)$. Then $F(n)=\bigcup_{y\in G(n)} R^y$ is an increasing sequence of finite subsets of $X$, hence for some $m$ and all $n\geq m$ we have $F(m)=F(n)$. Then $G(n)\subseteq \bigcup_{x\in F(m)} R_x$ for all $n$, hence all $G(n)$ are subsets of a fixed finite set. This implies that the sequence $G(n)$ eventually stabilizes. Since this sequence was arbitrary, $Y$ is CF.
\end{Proof}

\paragraph{} A function is \emph{finite-to-one} if the preimage of every element of its range is finite. The following is a notable consequence of Lemma~\ref{L.BackAndForth}.

\begin{Lemma} \label{L.finite-to-one} Suppose $f\colon X\to Y$ is finite-to-one and surjective.
	Then  $X$ is CF if and only if  $Y$ is CF. If $Y$ is DF so is $X$, but the converse does not hold in general.
\end{Lemma}

\begin{Proof} Take $R\subseteq X\times Y$ to be the graph of $f$. Then all vertical sections of~$R$ are singletons and all horizontal sections are finite.
	Thus   Lemma~\ref{L.BackAndForth} applies to imply that $X$ is CF if and only if  $Y$ is CF.

If there is an injection from $\bbN$ into $X$, then the composition of this function with $f$ has an infinite range, and from it one defines an injection from $\bbN$ into $Y$.

Finally, the Russell set (Definition~\ref{RussellSetDef}) is DF and it clearly has $\bbN$ as a two-to-one image.
\end{Proof}

\subsection{Examples and Constructions}

\paragraph{}
The first example of a DF set was due to {\sc A.\ Fraenkel} \cite{FraenkelGrundlagen}, who, assuming a certain extension ZFA of ZF is consistent (ZFA was later proved to be equiconsistent with ZF) constructed
a model of ZF containing an amorphous set $X$.  This $X$ is not only DF, but even
PF.

\paragraph{}\label{CohenModel}
The most famous Dedekind-finite set was constructed by {\sc P.\ Cohen} in his work on the
independence of the Axiom of Choice and the Continuum Hypothesis \cite{CohenSet}: a canonical infinite
Dedekind-finite subset $X$ of $\R$ in a certain model of ZF.  This $X$ is weakly
Dedekind-infinite (\ref{DedInfRSubProp}) and of course can be totally ordered, hence has quite
different properties from {\sc Fraenkel's} example (although {\sc Cohen's} construction was modeled
on {\sc Fraenkel's}). Also note that these sets exist in different universes: {\sc Fraenkel's} model is a model with atoms while {\sc Cohen's} model is a model of ZF.

\subsection{Rigid Sets}

Recall that a {\em permutation} of a set $X$ is a bijection from $X$ to $X$.

\paragraph{}
\begin{Definition}\label{RigidSetDef}\index{rigid set}\index{set!rigid}
An infinite set $X$ is {\em rigid} if every permutation of $X$ moves only finitely many elements (i.e.\ is
the identity on the complement of a finite subset).
\end{Definition}

\paragraph{}
Any set equipotent with a rigid set is rigid.  A subset of a rigid set is rigid (since a permutation of a
subset extends to a permutation of the whole set).  Since $\N$ is not rigid (it has many permutations
moving infinitely many elements, even all elements), a rigid set must be Dedekind-finite.  But ``most''
Dedekind-finite sets are not rigid (e.g.\ if $X$ is any infinite set, then $X\times Y$ is not rigid for any $Y$
with more than one element).

Rigid sets are interesting in the context of associated Hilbert spaces, since a permutation of $X$
naturally defines a unitary operator on $\ell^2(X)$; this unitary is a finite-rank perturbation of a scalar
if and only if the permutation moves only finitely many elements of $X$.

\paragraph{}
Being rigid seems to be similar to being amorphous.  But it turns out that the notions are
distinct.  Say an infinite set $X$ is {\em strongly rigid} if in any partition of $X$ into nonempty
subsets, all but finitely many cells are singletons, and {\em strongly amorphous} if it is
amorphous and strongly rigid.  A strongly rigid
set is rigid.  It can be shown that {\sc Fraenkel's} set is strongly amorphous, hence rigid, and that
{\sc Cohen's} generic subset of $\R$ is rigid, but it is not amorphous.  There are amorphous sets
containing infinitely many pairwise disjoint two-element sets (\S\ref{SuperRussell}), and such a set is not rigid (just
interchange each pair of socks).  And there are rigid sets which can be written as a countable
disjoint union of three-element sets in some models of ZF.

\subsection{Russell Sets and Russell Cardinals}\label{S.RussellSets}

\paragraph{}
\begin{Definition}\label{RussellSetDef}
A set $X$ is a {\em Russell set}\index{Russell set} if it is Dedekind-finite but can be written as a
countable union of pairwise disjoint two-element sets (called {\em pairs of Russell socks}).\index{Russell socks}
A cardinal $\kappa$ is a {\em Russell cardinal}\index{cardinal!Russell} if it is the cardinal of a Russell set.
\end{Definition}

\paragraph{}\label{RussellNotCF}
Note that a Russell set $X$ is DF but not CF, and is not amorphous.  There is a permutation
of $X$ interchanging each pair of socks, so $X$ is not rigid either.  If $\kappa$ is a Russell
cardinal, then $\kappa\pm1$ is not a Russell cardinal, but is DF and not CF.  In fact,
$\kappa\pm n$ is a Russell cardinal if and only if $n$ is even.  Russell socks are a graphic example
of a countable union of countable, even finite, sets which is not countable.

\paragraph{}\label{RussellDDF}
If $X$ is a Russell set, then $X$ is DDF. To prove this, let  $X=\bigsqcup_{n\in\N}X_n$ be the decomposition of $X$ into pairs. Suppose that $f\colon X\to X$ is a surjection. In order to prove that it is injective, let
\[
A=\{m\mid (\exists k\geq 1)(\exists n)f^k[X_n]\text{ is a proper subset of }X_m\}.
\]
We claim that $A$ is finite. Otherwise, for every $m\in A$ choose the lexicographically minimal pair $k(m), n(m)$ such that $f^{k(m)}[X_{n(m)}]$ is a proper subset of $X_m$. Let $x_m$ be the unique element of $X_m\setminus f^{k(m)}[X_{n(m)}]$. Then $\{x_m\mid m\in A\}$ is a countably infinite subset of $X$; contradiction.

For $m\in \bbN\setminus A$ let $g(m)$ be the minimal $n$ such that $f[X_n]=X_m$. Note that $g\colon (\bbN\setminus A)\to \bbN$ is an injection.
Suppose for a moment that for some $m_0\in \bbN\setminus A$ we have $g(m_0)\neq m_0$. Fix $x_{0}\in X_{m_0}$, and for $k\geq 1$ let $x_{k}$ be the unique element of $X_{g^k(m_0)}$ such that $f^k(x_{k})=x_{0}$. Since $g$ is an injection, the sequence $(m_k)$ is infinite and again we have a countably infinite subset $\{x_k\mid k\in \bbN\}$ of $X$; contradiction.

We have proven that the set $B=\{m\in \bbN\setminus A\mid g(m)=m\}$ is cofinite in $\bbN$, and $f[X_m]=X_m$ for all $m\in B$. Therefore $f$ is an injection on the cofinite subset $\bigcup_{m\in B}X_m$ of $X$. Since it is surjection, $f$ is an injection as required.

% First note that $\{n\in\N\mid\forall m\in\N, f[X_n]\not\subseteq X_m\}$ is finite, otherwise we can choose from each such $X_n$ the member which is mapped to a lower-indexed pair. Next, note that $\{m\in\N\mid\forall n\in\N,f[X_n]\subseteq X_m\to f[X_n]\neq X_m\}$ is also finite, as we may choose from such $X_m$ the element which is not in $f[X_n]$ for the least $n$ possible. So for cofinitely many of the pairs, $f$ maps one pair to the other, which means it must be injective on those. Suppose now that $f$ is not an injection, this must be witnessed at some point on the part where $f$ just moves one pair to the other, which means that $f$ moves, perhaps, multiple pairs to the same $X_n$, but it does so bijectively for each pair. Now pick the least $n$ for which $X_n$ is witnessing this phenomenon, let $a_0\in X_n$ be one of the points, and now simply observe that $a_n=f^{-n}(a_0)$ must be an injective sequence which defines a choice function from infinitely many pairs. Therefore $f$ must be injective and $X$ is a DDF set.

\paragraph{}
There are models of ZF containing Russell sets; for example, the second Cohen’s model (\cite[\S 5.4]{jech2008axiom}).  For a thorough and entertaining description of
Russell sets, including motivation for the sock terminology, see \cite{HerrlichTNumber}; some results
described there graphically show that the world of DF cardinals (when they exist at all) is
even more bizarre than we have indicated (for example, there is a family of $2^{\aleph_0}$
pairwise incomparable DF cardinals).  For one interesting observation, if one has a set of Russell
socks and adds another pair of socks,
the \emph{same} number ($\aleph_0$) of pairs of socks, but more socks (a strictly larger cardinality).  The latter fact is a consequence of Corollary~\ref{C.DF+finite}.

\paragraph{}
There are various models of ZF in which all the types of DF sets described above occur
(assuming there is a model for ZF at all, i.e.\ that ZF is consistent); see \S\ref{S.ModelsofZF}.  There is even a single model of ZF in which all such sets occur, including counterexamples for the false implications between the properties.

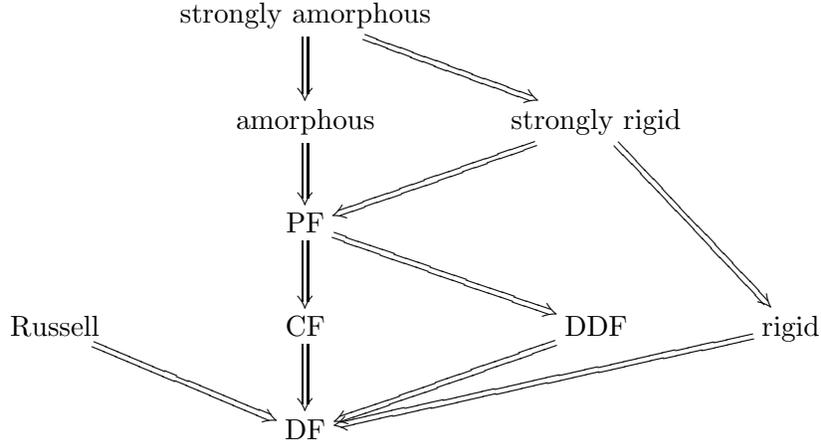
\begin{figure}[h]

	$$\xymatrix{ & \mbox{strongly amorphous}%\ar@{=>}[dl]
		\ar@{=>}[d]\ar@{=>}[dr] & & \\
		%\mbox{Hilbert-amorphous}\ar@{=>}[d]\ar@{=>}[dr] & & & \\
		 %\mbox{HPF}\ar@{=>}[d]\ar@{=>}[dr]
		 & \mbox{amorphous}\ar@{=>}[d] & \mbox{strongly rigid}\ar@{=>}[ddr]\ar@{=>}[dl] & \\
		%\mbox{HCF}\ar@{=>}[dr]
		& \mbox{PF}\ar@{=>}[d]\ar@{=>}[dr] & & \\ \mbox{Russell}\ar@{=>}[dr] & \mbox{CF}\ar@{=>}[d] & \mbox{DDF}\ar@{=>}[dl] & \mbox{rigid}\ar@{=>}[dll] \\
		& \mbox{DF} & & }$$
	\caption{\label{Fig.DF}Implications between different flavors of Dedekind--finite sets discussed in this section. See also Fig.~\ref{Fig.DF2} in \S\ref{S.HDF}.}
\end{figure}

\paragraph{}
Figure \eqref{Fig.DF} shows the flavors of DF sets considered in this paper.

\section{$\ell^2(X)$ For DF Sets}\label{S.ell2(X)}\setcounter{paragraph}{0}
\label{S.ell2}

If $X$ is a DF set, the Hilbert space $\ell^2(X)$ has some unusual properties depending on
the flavor of DF.

\smallskip

\paragraph{}
\begin{Definition}
	Let $X$ be a set, $\eta\in\ell^2(X)$, and $\eta(x)=\langle\eta,\xi_x\rangle$ the $x$'th coordinate of
	$\eta$ for each $x\in X$.  The {\em support} of $\eta$ is $supp(\eta)=\{x\in X:\eta(x)\neq0\}$.
\end{Definition}

For the following one needs to keep in mind that the union of a sequence of finite sets need not be countable (see Definition~\ref{RussellSetDef}).
\paragraph{}
\begin{Proposition}\label{VectSuppC}
	Let $X$ be a set, $\eta\in\ell^2(X)$, and $\eta_x=\langle\eta,\xi_x\rangle$ the $x$'th coordinate of
	$\eta$ for each $x\in X$.  Then $\supp(\eta)$ is the union of a sequence of finite subsets of $X$.
\end{Proposition}

\begin{Proof}
	Let $\epsilon>0$, and set $A_\epsilon=\{x:|\eta(x)|\geq\epsilon\}$.  If $B=\{x_1,\dots,x_n\}$ is an $n$-element subset
	of $A_\epsilon$, set
	$$\zeta=\sum_{k=1}^n \eta(x_k)\xi_{x_k}\ .$$
	then
	$$\|\eta\|^2\geq\|\zeta\|^2=\sum_{k=1}^n|\eta(x_k)|^2\geq n\epsilon^2$$
	so $n\leq\left ( \frac{\|\eta\|}{\epsilon}\right ) ^2$.  Thus $A_\epsilon$ is a finite set.  We have
	$\supp(\eta)=\cup_n A_{1/n}$.
\end{Proof}

Suppose $X$ is CF.  We then have severe restrictions on the support of vectors (note that (ii) of Proposition~\ref{StrongDFFiniteSupport} says that the canonical orthonormal basis of $\ell^2(X)$ is a Hamel basis; we will return to this observation in \ref{L.SDFBasis}).\footnote{Parts of Proposition~\ref{StrongDFFiniteSupport} have been proved in \cite{BrunnerLinear} in a slightly more general context.}

\paragraph{}
\begin{Proposition}\label{StrongDFFiniteSupport}
For every set $X$ the following are equivalent.
\begin{enumerate}
\item [(i)] $X$ is Cohen-finite.
\item [(ii)] Every vector in $\ell^2(X)$ has finite support.
\item [(iii)] Every sequence of disjointly supported unit vectors in $\ell^2(X)$ is finite.
\item [(iv)] Every sequence of
vectors in $\ell^2(X)$ has finite common support.
\item [(v)] The closed unit ball of $\ell^2(X)$ is sequentially compact.
\item[(vi)]  The standard orthonormal basis for $\ell^2(X)$ is a Hamel basis.

\end{enumerate}
\end{Proposition}

\begin{Proof} (i) implies (ii) follows immediately from \ref{VectSuppC} and \ref{StrongDFSubsets},
and (vi) is obviously equivalent to (ii).

(ii) implies (iii): If $\eta_n$ are disjointly supported unit vectors, then $\eta=\sum_{n=0}^\infty 2^{-n-1} \eta_n$ is a unit vector with infinite support.

(iii) implies (i): Suppose $X$ is not CF. Then there is a sequence of pairwise disjoint nonempty finite subsets of $X$, $F(n)$ for $n\in \N$. For each $n$ let $\eta_n=|F(n)|^{-1/2}\sum_{x\in F(n)} \xi_x$. These are unit vectors with disjoint supports.

It remains to prove (iv) and (v) are equivalent to the first three statements.  (iv) implies (ii)
is trivial.

(i) and (ii) together imply (iv):
If $(\eta_n)$ is a sequence of vectors in $\ell^2(X)$, then by (ii) $(\supp(\eta_n))$ is a sequence of finite subsets of $X$. Since $X$ is CF, this sequence can contain only finitely many distinct terms.

(iv) implies (v):
If $(\eta_n)$ is a sequence in $\ell^2(X)$, then the $\eta_n$ all lie in a finite-dimensional
subspace of $\ell^2(X)$.  But the closed unit ball of a finite-dimensional normed vector space
is compact (see \ref{FinDimHilbSp} and \ref{InfDimTotBdd}).

(v) implies (iii): Suppose that (iii) fails. An infinite sequence of unit vectors with disjoint supports does not have a norm-convergent subsequence (e.g.\ because it weakly converges to $0$).
\end{Proof}

\paragraph{}
Thus if $X$ is CF, the closed unit ball of $\ell^2(X)$ is a $\sigma$-complete metric
space which is sequentially compact but not compact (not totally bounded, cf.\ \ref{InfDimTotBdd}).

\paragraph{}
In fact, if $X$ is CF, $\ell^2(X)$ is the union of the finite-dimensional subspaces corresponding
to finite subsets of $X$.  There are two natural topologies on $\ell^2(X)$: the norm topology and the weak topology.\footnote{This is the pointwise convergence topology when vectors are identified with the associated linear functionals via the Riesz Representation Theorem, Theorem~\ref{RieszRepThm}.}  There is also the inductive limit topology from the finite-dimensional subspaces.
It may be that these three topologies coincide, at least on bounded subsets of $\ell^2(X)$.
(Note that Tikhonov's Theorem is needed to prove that the closed unit ball of a Hilbert space is weakly compact---an instance of the Banach--Alaoglu Theorem---so in this setting we do not necessarily have weak compactness of the ball.)

\bigskip

From this, we can obtain a nice example from \cite{BrunnerSB} of an inner product space which is Cauchy-complete but not $\sigma$-complete (another claimed example from \cite{BrunnerSB}
is wrong, cf.\ \ref{RussellOrthSeq}):

\paragraph{}
\begin{Example}\label{CauchyCompleteEx}
Let $X$ be a CF set.  Every vector in $\ell^2(X)$ is of the form $\displaystyle{\sum_{x\in X} c_x\xi_x}$
with the $c_x\in\bbC$ and only finitely many $c_x$ nonzero (Proposition~\ref{StrongDFFiniteSupport}).  Let
$$ \calH_0=\left \{ \sum_{x\in X}c_x\xi_x\in\ell^2(X):\sum_{x\in X}c_x=0\right \} \ .$$
Then $ \calH_0$ is a vector subspace of $\ell^2(X)$.  If $(\eta_n)$ is a Cauchy sequence in $ \calH_0$,
then there is a finite subset $Y$ of $X$ such that all the $\eta_n$ are supported in the span $\cY$
of $\{\xi_y:y\in Y\}$.  Since $\cY$ is finite-dimensional, $\eta_n\to\eta$ for some $\eta\in\cY$,
and $\eta\in \calH_0$.  Thus $ \calH_0$ is Cauchy-complete.  However, $ \calH_0$ is not closed in
$\ell^2(X)$, and is in fact dense.  For if $\displaystyle{\zeta=\sum_{x\in X} c_x\xi_x}$ is in
$ \calH_0^\perp$, all the $c_x$ must be equal (just compare with vectors in $ \calH_0$ which are
differences of two standard basis vectors), hence all $c_x$ must be 0, so $ \calH_0^\perp=\{0\}$.
Alternatively, define a linear functional $\phi$ on $\ell^2(X)$ by
$$\phi \left ( \sum_{x\in X}c_x\xi_x\right ) =\sum_{x\in X}c_x\ .$$
It can be shown directly that $\phi$ is unbounded, and hence $ \calH_0=\cN(\phi)$ (the kernel of $\phi$) is dense in
$\ell^2(X)$.  Thus $ \calH_0$ is not $\sigma$-complete.
\end{Example}

\subsection{$\ell^2(X)$ for Russell Sets}

\paragraph{}\label{RussellOrthSet}
We now let $X$ be a DF set which is not CF, i.e.\ $X$ contains a sequence $(A_n)$
of pairwise disjoint nonempty finite sets.  For each $n$ set
$$\eta_n=\frac{1}{\sqrt{|A_n|}}\sum_{x\in A_n}\xi_x\ .$$
Then $(\eta_n)$ is a sequence of unit vectors in $\ell^2(X)$ with disjoint supports, hence orthonormal.
Thus, although $\ell^2(X)$ has an orthonormal basis whose cardinality is DF,
$\ell^2(X)$ contains an orthonormal sequence of vectors (whose cardinality is $\aleph_0$, hence
not DF).  In fact, see \ref{NonuniqueBasis}.

\paragraph{}\label{RussellOrthSeq}
Now let $X$ be a Russell set, with $|A_n|=2$ for all $n$.  Then $\eta_n$ is the (normalized) sum of the basis
vectors in each pair of socks.  Let $\cY$ be the closed subspace of $\ell^2(X)$ spanned by the
$\eta_n$.  Then $\cY$ is isometrically isomorphic to $\ell^2(\N)$.

\smallskip

\noindent
Let $\cZ=\cY^\perp$.  Then $\cZ$ is the closed span of the differences of the vectors in each pair
of socks.  These differences are only well defined up to sign, and there does not exist a global
choice of signs; however, the one-dimensional subspace spanned by the difference is well defined
for each pair.  Thus $\cZ$ is spanned by a sequence of one-dimensional subspaces, although there
is not a global choice of unit vectors in these subspaces (since we are working with complex Hilbert spaces, we have a circle of unit vectors in each subspace, but splitting the circle into left and right
halves would give a global choice of signs).  But it seems
unlikely that $\cZ$ has an orthonormal basis; in fact, it seems likely that any orthonormal set in
$\cZ$ is finite (this is claimed in \cite{BrunnerSB}, with reference to a proof in \cite{BrunnerFraenkel}
which we are unable to follow or see its relevance; it is also claimed that this proof shows that $\cZ$ is not $\sigma$-complete, which is wrong since it is a closed subspace).  So $\cZ$ may be a Hilbert space in which no infinite-dimensional subspace has
an orthonormal basis. In \ref{L.ExistenceOfBasis} we will prove that with a stronger assumption on $X$ this is indeed the case. By Theorem~\ref{SepHilbSpaceThm}, neither $\cZ$ nor $\ell^2(X)$ is separable, although $\ell^2(X)$ is the ``direct sum'' of countably many 2-dimensional Hilbert spaces (infinite direct sums and products, and inductive limits, are not well defined without some Choice, cf.\ \ref{IndLim}).

It may be worth noting that if $\cZ$ has a basis $\cB$, then $\ell^2(X)$ has two bases of different cardinalities, namely the canonical basis and $\{\eta_n: n\in \N\}\cup \cB$ (compare with \ref{L.UniquenessOfBasis}). Moreover, regardless of whether $\cZ$ has a basis or not,~$\ell^2(X)$ has two bases of different cardinalities; see Example~\ref{NonuniqueBasis}.

\smallskip

\noindent
Note that there is a self-adjoint unitary operator on $\ell^2(X)$ which interchanges the basis
vectors in each pair of socks, for which $\cY$ and $\cZ$ are the $+1$ and $-1$ eigenspaces
respectively.

\subsection{Bases in Hilbert Spaces}

\paragraph{} The Axiom of Choice implies that every Hilbert space has an orthonormal basis, and that all of its bases have the same cardinality. We will see that the first statement can fail, and give  positive and negative results about the second.

%\paragraph{}\label{SuperRussell} For the following we need a strong assumption (stronger than  the existence of a Russell set, \ref{RussellSetDef}):
%There exists a set $X$ which can be presented as $X=\bigsqcup_{z\in Z} X_z$, where $Z$ is strongly amorphous, and $X_z$ are pairwise disjoint two-element sets, and there is no choice function $f\colon Z\to X$ (this is consistent with ZF by e.g. (REF:Truss)).  In this situation $Z$ is automatically CF. Since the function that collapses $X_z$ to $z$ is finite-to-one, \ref{L.finite-to-one} implies that $X$ itself is CF.

\paragraph{}\label{SuperRussell}  For the following we need a strong assumption (stronger than  the existence of a Russell set, \ref{RussellSetDef}):
\begin{quote}
	There exists an amorphous set $X$ which can be presented as $X=\bigsqcup_{z\in Z} X_z$, where $Z$ is strongly amorphous and $X_z$ are pairwise disjoint two-element sets.
\end{quote}
Relative consistency of the existence of such $X$ with ZF is given in \cite[Theorem~6.3]{truss1995structure}, but in Proposition~\ref{L.Truss} we show why this follows from general theory as presented in  \cite{karagila2022choiceless}.
If $X$ has  this property, then  there is no choice function $f\colon Z\to X$, because the range of such function would be an infinite, co-infinite subset of $X$. Moreover,  $X$ is amorphous. [If $Y$ is a subset of $X$, then for all but finitely many $z\in Z$, $Y$ either includes~$X_z$ or is disjoint from it. Since $Z$ is amorphous, the set $\{z\mid X_z\subseteq Y\}$ is finite or cofinite, thus $Y$ is cofinite or cofinite.]  Since $X$ is amorphous, it is also CF.  It is not rigid  because the permutation that interchanges the elements of $X_z$ for all $z$ has no fixed points.%Also, the index set  $Z$ is automatically CF. Since the function that collapses $X_z$ to~$z$ is finite-to-one, Lemma~\ref{L.finite-to-one} implies that $X$ itself is CF.

%\paragraph{}\begin{Theorem} \label{L.ExistenceOfBasis} Suppose $X$ is as in \ref{SuperRussell}. Then some (necessarily complemented, by Corollary~\ref{OrthDecompProp}) subspace of $\ell^2(X)$ has no basis.
% \end{Theorem}

\paragraph{}\begin{Theorem} \label{L.ExistenceOfBasis} Suppose $X$ is as in \ref{SuperRussell}. Then the following holds.
	\begin{enumerate}
		\item \label{1.L.existence.of.basis} $\ell^2(X)$ has a closed subspace $\cY$ isomorphic to $\ell^2(Z)$, and there is no injection from $Z$ into $X$.
		\item \label{2.L.existence.of.basis} The orthogonal complement of $\cY$   has no basis.
		\item \label{3.L.existence.of.basis}The orthogonal complement of $\cY$ contains no infinite orthonormal set.
	\end{enumerate}
\end{Theorem}

\begin{Proof}  \eqref{1.L.existence.of.basis}  Write $X=\bigsqcup_{z\in Z} X_z$.
	The canonical basis of $\ell^2(X)$ is  $\xi_{z,x}$, for $z\in Z$ and $x\in X_z$. (Clearly the index $z$ is redundant, but it will come handy in the proof.)  Consider the subspace $\cY$ of $\ell^2(X)$  spanned by the vectors $\eta_z=2^{-1/2}\sum_{x\in X_z}\xi_{z,x}$. Since each $X_z$ has exactly two elements, each~$\eta_z$ is a unit vector. These vectors have disjoint supports and are therefore orthogonal.  Clearly their closed linear span, $\cY$,  is isomorphic to $\ell^2(Z)$.

	Suppose that $f\colon Z\to X$ is an injection from $Z$ into $X$. Since $Z$ is amorphous, the set $\{z\mid f(z)\in X_z\}$ is either finite or cofinite. By the assumption that there is no choice function on $Z$, it is finite. Therefore the function $g\colon Z\to Z$ uniquely defined by the requirement $f(z)\in X_{g(z)}$ has at most finitely many fixed points. This function also has the property that the preimage of every element has at most two points; we will say that $g$ is \emph{at most two-to-one}.

	For $z\in Z$ consider the forward orbit $\{g^n(z)\mid n\geq 1\}$. It cannot be infinite, since otherwise there would be an injection from $\bbN$ into~$Z$. Therefore for all but finitely many  $z\in Z$ there are  $m(z)<n(z)$ such that $g^{m(z)}(z)=g^{n(z)}(z)$. Choose this pair so that $(m(z),n(z))$ is the minimal possible (in the lexicographic order).  Since $Z$ is amorphous, for every $S\subseteq \bbN^2$, the set $Z_S=\{z\in Z\mid (m(z),n(z))\in S\}$ is finite or cofinite. Therefore
	\[
	\{S\subseteq \bbN^2\mid Z_S\text{ is cofinite in $Z$}\}=\{S\subseteq \bbN^2\mid Z_S\text{ is infinite}\}
	\]
	and this set is a filter. Since $\bbN$ is not amorphous, it is a principal filter, and there is a pair $(m,n)$ such that
	\[
	Z'=\{z\in Z\mid (m(z),n(z))=(m,n)\}
	\]
	is cofinite in $Z$ (and therefore amorphous).

	We claim that $m=1$.
	For $k\geq 1$ let
	\[
	Z[k]=\{z\in Z’\mid (\exists y\in Z) g^k(y)=z\}.
	\]
	If $1\leq k$ then $g[Z[k]]=Z[k+1]$. Since $g$ is at most two-to-one, if $m>1$ and $Z[m]$ is infinite then $Z[m-1]$ is infinite as well. By the minimality of $m$ we conclude that $m=1$.

	Therefore $g(z)\in X_z$ for all but finitely many $z\in Z$. This contradicts the fact that $g$ has at most finitely many fixed points proven earlier.

\eqref{3.L.existence.of.basis} We claim that the orthogonal complement of $\cY$ in $\ell^2(X)$ includes no infinite orthonormal set.  Assume otherwise and fix one, denoted $\cB$.  Let
	\[
	R=\{ (z,\eta)\in Z\times\cB : \langle \eta, \xi_{z,x}\rangle \neq 0\text{ for some $x\in X_z$}  \}.
	\]
	We will use Lemma~\ref{L.BackAndForth} and the notation $R_z$, $R^\eta$ for horizontal and vertical sections of $R$ introduced there.
	Since $X$ is CF, every vector in $\ell^2(X)$ has finite support by Proposition~\ref{StrongDFFiniteSupport}. This means that  for every $\eta$ in $\cB$ the horizontal section $R^\eta$ of $R$ is finite.

	We claim that  for  every $z\in Z$, the vertical section $R_z$ of $R$ is finite.  To see this, fix $z$ and let $P_z$ denote the projection to the span of $\xi_{z,x}$, for  $x\in X_z$.  For every $n$, the set $\{\eta \in \cB: \|P_z(\eta)\|>1/n\}$ is finite. Since $X$ is CF, there is $n$ such that $P_z(\eta)$ is nonzero if and only if  $\|P_z(\eta)\|>1/n$. Since the set of $\eta$ such that $\|P_z(\eta)\|>1/n$ is finite,  we conclude that $R_z$ is finite for all $z\in Z$.

	Therefore the assumptions of Lemma~\ref{L.BackAndForth} (i) are satisfied. This lemma implies that there are a set~$W$ and partitions $Z=\bigsqcup_{w\in W} Z_w$ and $\cB=\bigsqcup_{w\in W} \cB_w$ into finite sets such that $\bigcup\{R_x: x\in Z_w\}\subseteq \cB_w$  and $\bigcup\{R^\eta: \eta\in \cB_w\}\subseteq Z_w$ for all $w\in W$.

	Then  $Z_w$, for $w\in W$, are equivalence classes of an equivalence relation on $Z$. Since $Z$ is strongly amorphous, all but finitely many of the classes are singletons. Let $Z'$ be the union of the singleton classes. Since $\ell^2(X_z)$ is two-dimensional, for each $z\in Z'$ there is a unique $\eta(z)\in\cB$   with $R^{\eta(z)}=\{z\}$. For each $z\in Z$ we have  $\eta(z)=\sum_{x\in X_z} \lambda_x \xi_{z,x}$. Since $\eta(z)$ is orthogonal to $2^{-1/2}\sum_{x\in X_z} \xi_{z,x}$ we have $\sum_{x\in Z_x} \lambda_x=0$.   We can therefore define a function $f\colon Z'\to X$ by
	\[
	f(z)=x\quad\text{if and only if}\quad 0\leq  \arg(\lambda_x)<\pi.
	\]
	Then $f$ is a choice function on $Z'$. Since $Z\setminus Z'$ is finite, it can be extended to a choice function on all of $Z$, contradicting our assumption.

	\eqref{2.L.existence.of.basis} clearly follows from
	\eqref{3.L.existence.of.basis}.
\end{Proof}

\paragraph{}
The space without a basis constructed in \ref{L.ExistenceOfBasis} is ``equal'' to $\bigoplus_{z\in Z} \calH_z$ (recall that infinite direct sums of Hilbert spaces are not really well defined; we just mean the union of the finite direct sums is dense), where (using the notation from the proof of \ref{L.ExistenceOfBasis}) for every $z$,  $ \calH_z=\ell^2(X_z)\ominus \Span(\sum_{x\in X_z}\xi_{z,x})$ is a one-dimensional Hilbert space.

\paragraph{}\begin{Proposition} \label{AlmostExistenceOfBasis} Suppose that $ \calH$ is a real Hilbert space and there is a family $\calH_z$, for $z\in Z$, of one-dimensional, pairwise orthogonal, subspaces of $ \calH$ such that $\bigcup_{F\in \Fin(Z)} \Span\bigcup\{\calH_z\mid z\in F\}$ is dense in $ \calH$.  Then $ \calH$ is a subspace of a Hilbert space with a basis.
\end{Proposition}

\begin{Proof} For each $z$, the unit ball $S(z)$ of $ \calH_z$ has two elements. Let $Y=\sqcup_{z\in Z} S(z)$. Then $v_z\colon \calH_z\to \ell^2(S(H_z))$ defined by $v_z(\xi)=2^{-1/2} \sum_{\eta\in S(z)} \langle \xi, \eta\rangle \eta$ is an isometry for every $z$, and $\bigoplus_{z\in Z}v_z\colon H\to \ell^2(Y)$ is an isometry.
\end{Proof}

\paragraph{} It is not clear how to extend the argument of \ref{AlmostExistenceOfBasis} to any other space of the form $\bigoplus_{z\in Z} \calH_z$, even if $\dim( \calH_z)=2$ for all $z$ or if $\dim( \calH_z)=1$ but it is a complex vector space. (The latter is  the case of \ref{L.ExistenceOfBasis}, but the space constructed in this proposition is rather specific.)

Compare the following with Proposition~\ref{StrongDFFiniteSupport} of which it is a continuation. Parts of it appear in \cite{BrunnerLinear} in a slightly more general context.

\paragraph{} \label{L.SDFBasis}
\begin{Proposition} For every infinite-dimensional Hilbert space $ \calH$ with an orthonormal basis the following are equivalent.
\begin{enumerate}
\item [(i)] $ \calH$ has a CF orthonormal basis.
\item [(ii)] Every infinite orthonormal set in $ \calH$ is CF.
\item [(iii)] Every orthonormal basis is a Hamel basis.
\end{enumerate}
\end{Proposition}

\begin{Proof} (ii) trivially implies (i).

(i) implies (ii): If (i) holds, then by Proposition~\ref{OrthBasisL2} we may assume $ \calH=\ell^2(X)$ for a CF set~$X$. Fix an orthomormal set $\cB$ in $ \calH$. Let
\[
R=\{(x,\eta)\in X\times \cB: \langle \xi_x,\eta\rangle \neq 0\}.
\]
For every $\eta\in \cB$ the horizontal section $R^\eta$ is equal to the support of $\eta$, and therefore finite by Proposition~\ref{StrongDFFiniteSupport}.

Fix $x\in X$. We will prove that $R_x$ is finite. Assume otherwise. Since $\sum_{\eta\in R_x} |\langle \xi_x,\eta\rangle|^2=1$, the set $G(n)=\{\eta\in \cB: |\langle \xi_x,\eta\rangle|\geq 1/n\}$ is finite for all $n$. Therefore $F(n)=\bigcup_{\eta\in G(n)} R^\eta$ is a nondecreasing sequence of finite subsets of $X$. Since $X$ is CF, there exists $m$ such that $F(n)=F(m)$ for all $n\geq m$.

Every $G\subseteq \cB$ is included in the linear span of $\{\xi_x: x\in \bigcup_{\eta\in G}R^\eta\}$. This implies that $|G(n)|\leq |F(n)|$ (conveniently, this is the additional assumption of Lemma~\ref{L.BackAndForth} (ii)) which is by the previous paragraph $\leq |F(m)|$ for all $n$. Since the sequence $G(n)$ is increasing, it stabilizes. Hence its union is finite, contradicting the assumption that $R_x$ was infinite.

Since $x\in X$ was arbitrary, all horizontal sections of $R$ are finite, and   Lemma~\ref{L.BackAndForth} (ii) implies that $\cB$ is CF.

(i) implies (iii): Assume (i) and fix an orthonormal basis of $ \calH$. By (ii) we may assume $ \calH=\ell^2(X)$ for a CF set $X$. Proposition~\ref{StrongDFFiniteSupport} implies that every vector in $\ell^2(X)$ is a linear combination of finitely many basis vectors.

(iii) implies (i): Assume (i) fails, hence $ \calH=\ell^2(X)$ and $X$ is not CF. By Proposition~\ref{StrongDFFiniteSupport} (ii) implies (i), some vector in $ \calH$ has infinite support, thus $\xi_x$, for $x\in X$, is not a Hamel basis.
\end{Proof}

%\subsection{Orthogonal Dimension}
%
%\paragraph{}
%There remains the rather subtle question of whether $\ell^2(X)\cong\ell^2(Y)$ implies $|X|=|Y|$,
%i.e.\ whether the cardinality of an orthonormal basis (when there is one!) for a Hilbert space is well defined.  It is under
%the AC, but without Choice the answer is no (\ref{NonuniqueBasis}); how badly this fails is somewhat up in the air.
%There are some conditions under which uniqueness of orthogonal dimension holds (\ref{L.UniquenessOfBasis}).

\subsection{Uniqueness of Orthogonal Dimension}

Can we have a Hilbert space with bases of different cardinalities?  The answer is yes:

\paragraph{}
\begin{Example}\label{NonuniqueBasis}
Let $X$ be a Dedekind-finite set which is not Cohen-finite, and let $\cY$ be the closed subspace of $\ell^2(X)$ spanned by
an orthonormal sequence (\ref{RussellOrthSet}) and $\cZ=\cY^\perp$.  Then $\cY\cong\ell^2(\N)$.  If $X\sqcup\N$ is the disjoint union
of $X$ and $\N$, then
\begin{multline*}
\ell^2(X\sqcup\N)\cong\ell^2(\N)\oplus\ell^2(X)\cong\ell^2(\N)\oplus\cY\oplus\cZ\\
\cong\ell^2(\N)\oplus\ell^2(\N)\oplus\cZ\cong\ell^2(\N)\oplus\cZ\cong\cY\oplus\cZ\cong\ell^2(X)
\end{multline*}
since $\ell^2(\N)\oplus\ell^2(\N)\cong\ell^2(\N)$.  Thus, if $X$ is DF, $\ell^2(X)$ has an orthonormal basis indexed
by the DF set $X$, and an orthonormal basis indexed by the Dedekind-infinite set $X\sqcup\N$.

It follows that if $Y$ is $X$ with any finite number of points added or removed, i.e.\
$|Y|=|X|\pm n$ for some $n\in\N$, then $\ell^2(Y)\cong\ell^2(X)$.  For any $n\neq 0$, $|Y|\neq|X|$.
If $X$ is a Russell set, then $Y$ is a Russell set if and only if $n$ is even (\S\ref{RussellNotCF}).
\end{Example}

The following shows that the cardinalities of orthonormal bases cannot be too wild.

\paragraph{} \label{L.UniquenessOfBasis}
\begin{Theorem} Suppose that $\ell^2(X)$ and $\ell^2(Y)$ are isomorphic, $U\colon \ell^2(X)\to \ell^2(Y)$  is  a unitary,  and $X$ is CF.
\begin{enumerate}
\item [(i)]	 There are a set $Z$ and partitions $X=\bigsqcup_{z\in Z} X_z$ and $Y=\bigsqcup_{z\in Z} Y_z$ into finite sets such that
$U[\ell^2(X_z)]=\ell^2(Y_z)$ (hence $|X_z|=|Y_z|$)  for all $z\in Z$.
\item [(ii)] If $X$ is in addition strongly amorphous, then $|X|=|Y|$  and there are cofinite $X'\subseteq X$, $Y'\subseteq Y$, and a bijection $f\colon X'\to Y'$ such that $U\xi_x=c \xi_{f(x)}$ for some $c\in \bbC$ with $|c|=1$,  for all $x\in X'$.
\end{enumerate}
\end{Theorem}

\begin{proof} (i) By \ref{L.SDFBasis}, $Y$ is CF, hence   \ref{StrongDFFiniteSupport} implies that every vector in $\ell^2(X)$ and $\ell^2(Y)$ has finite support. Therefore all vertical and horizontal sections of $R\subseteq X\times Y$ defined by
\[
R=\{(x,y)\in X\times Y:  \langle U\xi_x, \xi_y\rangle\neq 0\}
\]
are finite. Since $X$ is CF, Lemma~\ref{L.BackAndForth} implies that there are a set $Z$ and partitions $X=\bigsqcup_{z\in Z} X_z$ and $Y=\bigsqcup_{z\in Z} Y_z$ into finite sets such that $\bigcup\{R_x: x\in X_z\}\subseteq Y_z$  and $\bigcup\{R^y: y\in Y_z\}\subseteq X_z$ for all $z\in Z$.

 Then the restriction of $U$ to $\ell^2(X_z)$ is an isometry onto $\ell^2(Y_z)$ for all $z\in Z$. Since these are finite-dimensional spaces, $|X_z|=|Y_z|$.
%
%For a moment fix $x\in X$. Recursively define $F(n)\subseteq X$ and $G(n)\subseteq Y$ as follows. Let $F(0)=\{x\}$, and for all $n$ let
%\begin{align*}
%G(n)&=\{y\in Y:  \langle U\xi_x, \xi_y\rangle\neq 0\text{ for some $x\in F(n)$}\},\\
%F(n+1)&=\{x\in X: \langle U\xi_x, \xi_y\rangle\neq 0\text{ for some $y\in G(n)$}\}.
%\end{align*}
%	Since $X$ and $Y$ are CF,  Lemma~\ref{StrongDFFiniteSupport}, every vector in $\ell^2(X)$ and $\ell^2(Y)$ has finite support. Therefore all $F(n)$ and all $G(n)$ are finite.
%	Clearly $F(n)\subseteq F(n+1)$ for all $n$. Since $X$ is CF, there exists $m=m(x)$ such that $F(m)=F(n)$ for all $n\geq m$.  This implies $G(m)=G(n)$ for all $n\geq m$. Let $F[x]=F(m)$ and $G[x]=G(m)$. Then the restriction of $U$ to $\ell^2(F[x])$ is an isometry onto $\ell^2(G[x])$. Since these are finite-dimensional spaces, $|F[x]|=|G[x]|$.

%	Fix $x'\in F[x]$. By induction on the minimal $n$ such that $x'\in F(n)$, one proves that $F_{x'}=F[x]$, and therefore $G_{x'}=G[x]$. Let $Z=\{F[x]: x\in X\}$. For every $z\in Z$ there is a unique finite subset $G(z)$ of $Y$ which satisfies $G=G[x]$ for some (every) $x$ such that $F[x]=z$.  Then $X_z=z$ and $Y_z=G(z)$ are the required partitions of $X$ and $Y$ into finite sets.

	To prove (ii), use the fact that since $X$ is strongly amorphous, $X_z$ is a singleton for a cofinite $Z'\subseteq Z$. Thus $Y_z$ is also a singleton for all $z\in Z'$, and we have a bijection between cofinite subsets of $X$ and $Y$. Using  $|X_z|=|Y_z|$ for all $z$ again, we  extend it to a bijection between $X$ and $Y$. Since $X$ is amorphous, for all but finitely many $x\in X$ we have $U\xi_x=c\xi_{f(x)}$ for some fixed $c\in\bbC$, $|c|=1$.
	\end{proof}

By applying (ii) of Theorem~\ref{L.UniquenessOfBasis} in the case when $X=Y$ is strongly amorphous, we obtain the following.

\paragraph{}\label{StronglyDFBasis} \begin{Corollary} If $X$ is strongly amorphous then every orthonormal basis of $\ell^2(X)$ includes $\{c\xi_x:x\in Y\}$ for some cofinite subset $Y$ of $X$ and some $c\in\bbC$, $|c|=1$.
\end{Corollary}

\subsection{Universality of  Hilbert Spaces With a Basis?}\label{S.universality}

\paragraph{}Can one prove in ZF that  every Hilbert space is the quotient of a Hilbert space with a basis, and/or that every Hilbert space is isomorphic to a subspace of a Hilbert space with a basis? We will see that these two questions are related (\ref{ACCQuotientSubspace}).  The only example of a Hilbert space without a basis known to us, given in \ref{L.ExistenceOfBasis},  is both a subspace and a quotient of a Hilbert space  with a basis (the latter follows by \ref{ACCQuotientSubspace}).

\paragraph{} A standard argument shows the following.

\begin{Lemma}\label{ACCquotient}
Every Banach space is the quotient of a Banach space with a basis.
\end{Lemma}

\begin{Proof}  Suppose $\cX$ is a Banach space and let  $Y$ be the unit sphere of $\cX$.  Denoting the canonical basis of $\ell^1(Y)$ by $\delta_y$, for $y\in Y$ define $T\colon \ell^1(Y)\to \cX$ by
\[
T(\sum_y \lambda_y \delta_y)=\sum_y \lambda_y y.
\]
To see that $T$ is well-defined, note that  every $\sum_y\lambda_y \delta_y$ satisfies
 $\|\sum_y \lambda_y \delta_y\|=\sum_y |\lambda_y|\geq \sum_y |\lambda_y| \|y\|\geq \|\sum_y y\|$. Therefore , we have  $\|T\|\leq K$.

By the choice of $Y$, $T$ is surjective.
\end{Proof}

\paragraph{} There is a lot of room for improvement in  Lemma~\ref{ACCquotient}. %Is it true for every $\sigma$-complete $X$, or even every Banach space $X$?
If $X$ is a Hilbert space, can we choose $Y$ to be a Hilbert space? If the set $Y$ in the proof of  Lemma~\ref{ACCquotient}  is a dense subset of the unit sphere of $X$, is the operator $T$ still surjective (the standard proof of this uses the axiom of Dependent Choice)?

In order to prove \ref{ACCQuotientSubspace} below, we need the following fact interesting in its own right.

\paragraph{} \begin{Proposition} \label{ACCShortExact} Every short exact sequence  of Hilbert spaces in which the connecting maps are isometries splits.
\end{Proposition}

\begin{Proof} Suppose that $0\to \cX\to  \calH\to \cY\to 0$ is a short exact sequence of Hilbert spaces with isometries as connecting maps.  Thus the image of $\cX$ is a closed subspace of $ \calH$. By Corollary~\ref{OrthDecompProp}, this subspace has an orthogonal complement $\cX^\perp$ and there is an orthogonal projection $P_{\cX}$ from $ \calH$ onto $\cX^\perp$. The restriction of the quotient map $q\colon  \calH\to \cY$ to $\cX^\perp$ is a surjective isometry,  and therefore its inverse  is continuous.
\end{Proof}

It would be of interest to generalize \ref{ACCShortExact} to the case when the connecting maps are arbitrary bounded linear operators. The Polar Decomposition Theorem (\ref{PolarDecomposition}) may be relevant.

The following is an immediate consequence of \ref{ACCShortExact}.

\paragraph{}\begin{Corollary} \label{ACCQuotientSubspace} Every Hilbert space which is the quotient of a Hilbert space with a basis is a subspace of a Hilbert space with a basis. \qed
\end{Corollary}

\section{Bounded Operators}\label{S.Bounded}

If $ \calH$ is a Hilbert space, denote by $\cB( \calH)$ the set of bounded operators on~$ \calH$.

\paragraph{}
\begin{Proposition}\label{BHsigmacomp}
If $ \calH$ is a Hilbert space, then $\cB( \calH)$ is $\sigma$-complete (i.e.\ a Banach algebra) under the operator norm.
\end{Proposition}

\begin{Proof}
Let $(A_n)$ be a decreasing sequence of closed bounded sets in $\cB( \calH)$ whose diameters
go to 0.  Fix a vector $\xi\in \calH$, and set $A_n^\xi=\{T\xi:T\in A_n\}$.  Then the sets
$(\overline{A_n^\xi})$ form a decreasing sequence of closed bounded subsets of $ \calH$ whose
diameters go to 0, so $\cap_n\overline{A_n^\xi}$ is a singleton we call $S\xi$.  This defines a
function $S: \calH\to \calH$, which is clearly linear.  If $\|T\|\leq M$ for all $T\in A_1$, then for any $\xi$
$\|T\xi\|\leq M\|\xi\|$ for all $T\in A_1$, and hence $\|S\xi\|\leq M\|\xi\|$, so $S$ is bounded and
$\|S\|\leq M$.  To show $S\in\cap_nA_n$,  let $\epsilon>0$, and fix $m$ with $\diam(A_m)\leq\epsilon$.
Let $n\geq m$.  If $\xi$ is any unit vector in $ \calH$, the set $\overline{A_n^\xi}$ has diameter $\leq\epsilon$,
hence $\|T\xi-S\xi\|\leq\epsilon$ for all $T\in A_n$.  Thus $\|T-S\|\leq\epsilon$ for every $T\in A_n$.
Since the $A_n$ are decreasing and closed, $S\in\cap_nA_n$.
\end{Proof}

\smallskip

This result holds more generally if $ \calH$ is just a ($\sigma$-complete) Banach space;
in fact, if $\cX$ and $\cY$ are normed vector spaces with $\cY$ $\sigma$-complete, then
$\cB(\cX,\cY)$ is $\sigma$-complete.

\bigskip

The next important fact is a corollary of the Riesz Representation Theorem (\ref{RieszRepThm}); the standard proof (cf.\ \cite[XVI.9.11.2]{BlackadarReal}) works
verbatim.

\paragraph{}
\begin{Corollary}\label{AdjointOp}
Let $ \calH$ be a Hilbert space, and $T\in\cB( \calH)$.  Then there is a unique $S\in\cB( \calH)$ with
$$\langle Tx,y\rangle=\langle x,Sy\rangle$$
for each $x,y\in \calH$.  Write $T^*$ for this $S$.
\end{Corollary}

\paragraph{}\label{BHCAlg}
The adjoint operation $T\mapsto T^*$ has the usual properties (\cite[I.2.3.1]{BlackadarOperator}), which show that
$\cB( \calH)$ is a C*-algebra.  \ref{AdjointOp} is what allows
reasonable operator theory to be done on $ \calH$.

\paragraph{}
In the example of a Cauchy-complete, not $\sigma$-complete space $ \calH_0$  obtained from a CF set $X$ in \ref{CauchyCompleteEx}, all the conclusions \ref{HilbCloseApproxThm}, \ref{OrthDecompProp},
and \ref{RieszRepThm} fail.  Operators on $ \calH_0$ do not have
adjoints in general, or more precisely the adjoint of a bounded operator on $ \calH_0$ is a (bounded) operator on $ \calH$, the $\sigma$-completion of $ \calH_0$.
Here is a simple example of an operator on $ \calH_0$ with no adjoint.  Let $\zeta$ be a unit vector in $\ell^2(X)\setminus  \calH_0$, $\eta$ a unit vector in $ \calH_0$, and let $T(\xi)=\langle \xi,\zeta\rangle\eta$ for $\xi\in \calH_0$.  This is a rank one operator on $ \calH_0$, but its adjoint (in $\cB(\ell^2(X))$) is a rank-one operator whose range is the span of $\zeta$, so it is not an operator on $ \calH_0$.
Note that this example also gives failure of the Riesz Representation theorem for~$ \calH_0$ and  a bounded sesquilinear form on $ \calH_0$ that does not correspond to a bounded linear operator on $ \calH_0$.

\paragraph{}
As usual, we define a {\em positive operator} on $ \calH$ to be a $T\in\cB( \calH)$ with $T=T^*$
and $\langle T\xi,\xi\rangle\geq0$ for all $\xi\in \calH$.  If $T\in\cB( \calH)$, then $T^*T\geq0$.  The
converse also holds: in fact, every $T\geq0$ has a unique positive square root, an $S\geq0$
with $S^2=T$ (write $S=T^{1/2}$).  This is a special case of continuous functional calculus of self-adjoint elements. %In a paper in preparation  (\cite{BlackadarF}) we will show that a standard application of Shoenfield's Absoluteness Theorem implies that this  holds without any Choice, but
There is a direct proof not requiring any
Choice, cf.\ \cite{RieszN} or \cite{HalmosHilbert}.

%\paragraph{}\label{ContinuousFunctionalCalculus} ????? It is worth mentioning that the continuous functional calculus for normal operators on a Hilbert space (and even for commuting finite sets of normal operators on a Hilbert space) can be developed in ZF. First, the \cstar-algebra $C(K)$ of all continuous complex valued functions on a compact Hausdorff space $K$ is $\sigma$-complete with respect to the sup norm. Second, if $K$ is a nonempty compact subset of $\bbC^n$, then the $^*$-polynomials in $n$ commuting variables form a dense subset of $C(K)$. Thus for every $n$-tuple of commuting normal operators $\bar T$ with joint spectrum $K$ and every $f\in C(K)$ one can define $f(\bar T)$ as a uniform limit of polynomials.

\paragraph{} \begin{Theorem}{\sc [Polar Decomposition Theorem]}
\label{PolarDecomposition}
Let $ \calH$ be a Hilbert space, and $T$ a bounded linear operator on $ \calH$. Then there are a unique positive operator $|T|$ and a partial isometry $V$ such that $T=V|T|$. If we require $V$ to have the property that the range of $V^*V$ is $\cN(T)^\perp$, then $V$ is unique.
\end{Theorem}

\begin{Proof} Set $|T|=(T^*T)^{1/2}$. For every vector $\xi$, by a standard argument we have
\[
\||T|\xi\|^2=\langle |T|\xi, |T|\xi\rangle=
\langle T^*T\xi, \xi\rangle=\langle T\xi,T\xi\rangle=\|T\xi\|^2.
\]
Thus $|T|\xi\mapsto T\xi$ is an isometry from the range of $|T|$ to the range of $T$, and hence by
\ref{DenseIsom} extends to an isometry between the closures.  Extend this to $V\in\cB( \calH)$ by setting $V=0$ on $\Ran(|T|)^\perp=\cN(T)$ (the null space of $T$).   Then $\cN(V)=\cN(|T|)=\cN(T)$, and $T\xi=V|T|\xi$ for all $\xi\in  \calH$, as required.
\end{Proof}

\paragraph{}
As a consequence, many other standard properties of $\cB( \calH)$ persist in this setting, such as:
\begin{enumerate}
\item[(i)]  Every $T\in\cB( \calH)$ has left and right support projections.
\item[(ii)]  The projections in $\cB( \calH)$ form a complete lattice.
\item[(iii)]  The right annihilator of any subset of $\cB( \calH)$ is generated by a projection (i.e.\ $\cB( \calH)$ is an \AWstar-algebra).
\item[(iv)]  $\cB( \calH)$ admits not only continuous functional calculus of self-adjoint elements
\cite{BlackadarF}, but even Borel functional calculus; in particular, if $T=T^*\in\cB( \calH)$
and $B$ is a Borel subset of $\sigma(T)$, the spectral projection $E_B(T)$ is well defined with
the usual properties.
\item[(v)]   A bounded increasing net of positive operators converges in the strong operator
topology.
\end{enumerate}
 These results will appear in \cite{BlackadarF}.

\subsection{Compact Operators}

\paragraph{}
Under usual versions of Choice, there are several equivalent conditions which characterize
compact operators.  These conditions may not be equivalent in the absence of CC,
however, so the proper definition must be nailed down.  The following definition seems to be the
best and most useful one (whether or not the AC is assumed):

\paragraph{}
\begin{Definition}\label{CompactOpDef}
Let $ \calH$ be a Hilbert space, and $T: \calH\to \calH$ a linear operator.  Then $T$ is {\em compact}
if, whenever $A$ is a bounded subset of $ \calH$, $T(A)$ is totally bounded.  Denote by $\cK( \calH)$
the set of compact operators on $ \calH$.
\end{Definition}

\paragraph{}
If $T$ is a compact operator, the image of the unit ball of $ \calH$ is totally bounded, hence bounded,
so $T$ is bounded, i.e.\ $\cK( \calH)\subseteq\cB( \calH)$.  Since bounded sets in a finite-dimensional
normed vector space are totally bounded (Choice is not needed to prove this), every finite-rank
bounded operator is compact.  Using our definition, we have the next standard result (it is unclear
whether it holds using other definitions):

\paragraph{}
\begin{Proposition}\label{CompactOpIdeal}
Let $ \calH$ be a Hilbert space.  Then $\cK( \calH)$ is a two-sided ideal in $\cB( \calH)$, and is the
norm-closure of the set of finite-rank operators on~$ \calH$.
\end{Proposition}

\begin{Proof}
Since bounded operators send bounded sets to bounded sets and totally bounded sets to totally
bounded sets, $\cK( \calH)$ is closed under left and right multiplication by bounded operators.
$\cK( \calH)$ is obviously closed under scalar multiplication.  If $S$ and $T$ are compact operators,
$A$ is a bounded subset of $ \calH$, and $\epsilon>0$, cover $S(A)$ by open balls of radius
$\frac{\epsilon}{2}$ centered at $\{x_1,\dots,x_n\}$ and cover $T(A)$ by open balls of radius
$\frac{\epsilon}{2}$ centered at $\{y_1,\dots,y_m\}$; then the open balls of radius $\epsilon$
centered at the $x_j+y_k$ cover $S(A)+T(A)\supseteq(S+T)(A)$.  So $S+T$ is compact.  Thus $\cK( \calH)$
is a two-sided ideal in $\cB( \calH)$.

\smallskip

\noindent
Let $T$ be in the norm-closure of $\cK( \calH)$, $A$ a bounded subset of $ \calH$, and $\epsilon>0$.
Choose $S\in\cK( \calH)$ with $\|S-T\|<\frac{\epsilon}{2}$.  Cover $S(A)$ with open balls of radius
$\frac{\epsilon}{2}$ centered at $\{x_1,\dots,x_n\}$.  Then the balls of radius $\epsilon$ centered
at $\{x_1,\dots,x_n\}$ cover $T(A)$.  Thus $T\in\cK( \calH)$.  So $\cK( \calH)$ is norm-closed.

\smallskip

\noindent
Finally, let $T\in\cK( \calH)$ and $\epsilon>0$.  Cover the image of the closed unit ball of $ \calH$
under $T$ by open balls of radius $\frac{\epsilon}{2}$ centered at $\{x_1,\dots,x_n\}$.  Let $P$ be the
orthogonal projection of $ \calH$ onto the finite-dimensional subspace spanned by $\{x_1,\dots,x_n\}$.
[Existence of $P$ follows from \ref{OrthDecompProp}, but there is a more elementary argument:
any finite-dimensional subspace of an inner product space has an orthogonal complement by
applying the Gram-Schmidt process to an orthonormal basis of the subspace and a given vector in the space.]  Since $P$ has norm 1, it follows that $\|PT-T\|<\epsilon$, and
$PT$ has finite rank.  Thus $T$ is in the closure of the finite-rank operators.
\end{Proof}

\paragraph{}
\begin{Corollary}\label{CompactOpAdjoint}
Let $ \calH$ be a Hilbert space, and $T\in\cK( \calH)$.  Then $T^*$ is also in $\cK( \calH)$.
\end{Corollary}

%\begin{Proof}
%By density and norm-continuity of $T\mapsto T^*$, it suffices to show this for $T$ finite-rank.  Let $\cY$ be the range of $T$.  Then
%$\cY^\perp$ is the null space of $T^*$, so the range of $T^*$ is the same as the range of $T^*|_{\cY}$,
%which is finite-dimensional.
%\end{Proof}

\paragraph{}
Thus $\cK( \calH)$ is a C*-subalgebra of the C*-algebra $\cB( \calH)$.  (Actually any closed two-sided
ideal in a C*-algebra is a C*-subalgebra, but this is an elementary observation in the case of $\cK( \calH)$.)

A \cstar-algebra is AF (\emph{approximately finite}) if it is an ``inductive limit''  of a directed system of finite-dimensional \cstar-algebras (since inductive limits are tricky and not well defined by Bratteli diagrams without the AC, we will
take this to mean there is a system of finite-dimensional C*-subalgebras, directed by inclusion, whose
union is dense).
\paragraph{}
\begin{Proposition}\label{P.K(H)simple}
If $ \calH$ is a Hilbert space, then $\cK( \calH)$ is a simple C*-algebra, which is nonunital if $ \calH$ is infinite-dimensional.
For every set $X$, $\cK(\ell^2(X))$ is an AF algebra.
\end{Proposition}

%\begin{Proof}
%	Only the last part may require a proof. For every finite $F\subseteq X$ we have that $\cB(\ell^2(F))$ is isomorphic to $M_n(\bbC)$, where $n=|F|$, and  naturally identified to a subalgebra of $\cK(\ell^2(X))$, and $\cK(\ell^2(X))$ is the inductive limit of these algebras.
%\end{Proof}

\paragraph{}
We may more generally define adjoints and compactness for operators between two Hilbert spaces
in the same way, with analogous properties.

\paragraph{}\label{IndLim}
If $X$ is a set of Russell socks, then $X$ is the union of an increasing sequence of finite sets,
and hence there is an increasing sequence of finite-dimensional matrix algebras, each embedded
as a corner in the next, whose union is dense in $\cK(\ell^2(X))$.
%(i.e.\ $\cK(\ell^2(X))$ is an``AF algebra'').
 $\cK(\ell^2(\N))$ has a
seemingly identical such sequence.  But $\cK(\ell^2(X))$ is not isomorphic to $\cK(\ell^2(\N))$;
in fact, $\cK(\ell^2(X))$ is not separable since no countable set can approximate the rank one
projections onto the spans of the standard basis elements (cf.\ the proof of \ref{SepHilbSpaceThm}).  Thus inductive limits
of sequences of C*-algebras, even finite-dimensional ones, are not well defined in the usual
sense without CC.

\paragraph{}There is a notable analogy of \ref{IndLim} to a result of Katsura (\cite[Theorem~6]{katsura2006non}) showing that nonisomorphic non-separable AF-algebras may have the same Bratteli diagram.

\section{$\cB(\ell^2(X))$ For DF Sets}\setcounter{paragraph}{0}\label{S.B(H)}

\paragraph{}
If $X$ is a set, then the Hilbert space $\ell^2(X)$ and the C*-algebras $\cB(\ell^2(X))$,
$\cK(\ell^2(X))$, and the Calkin algebra $\cQ(\ell^2(X))=\cB(\ell^2(X))/\cK(\ell^2(X))$ are invariants of the cardinal $|X|$, along with associated
objects like the lattice of ideals of $\cB(\ell^2(X))$.  Properties of these objects reflect properties
of $|X|$, in a more interesting and varied way than in the AC case.

\bigskip

We now examine properties of $\cB( \calH)$, $\cK( \calH)$, and the Calkin algebra $\cQ( \calH)$
for $ \calH=\ell^2(X)$, $X$ DF.

 A prominent open problem in theory of  \cstar-algebras is whether every stably finite \cstar-algebra has a tracial state; we recall the definitions for the benefit of the readers with stronger background in set theory than in operator algebras; note the analogy with the definition of Dedekind—finiteness.

\paragraph{} \begin{Definition}
Suppose that $A$ is a unital \cstar-algebra.  It is called \emph{finite} if there is no $v\in A$ such that $1_A=v^*v$ and $vv^*\neq v^*v$,\footnote{Such $v$, if it exists, is called a \emph{proper isometry} since its image under any GNS representation is a proper isometry.} \emph{infinite} if it is not finite, and \emph{stably finite} if $M_n(A)$ is finite for all $n$.
It is \emph{properly infinite} if there are no $v$ and $w$ in $A$ such that $v^*v=w^*w=1_A$ and $v^*w=0$.
\end{Definition}

See \cite[V.2.1]{BlackadarOperator} for more information.

One might hope (and we {\em did} hope) that if $ \calH=\ell^2(X)$ for a DF set $X$ then  $\cB( \calH)$ would be finite, i.e.\ have no nonunitary
isometries (even stably finite since
$ \calH^n\cong\ell^2(Y)$ with $Y=X\times\{0,\dots,n-1\}$ which is also DF).   But things are not so simple:

\paragraph{}
\begin{Example}\label{RussellEx}
Let $X$ be a set which is DF but not CF.  Then $\ell^2(X)\cong\cY\oplus\cY^\perp$ (\ref{RussellOrthSet}), where
$\cY\cong\ell^2(\N)$.  The unilateral shift $S$ is a nonunitary isometry on $\ell^2(\N)$, so
$S\oplus I$ is a nonunitary isometry in $\cB(\cY\oplus\cY^\perp)=\cB(\ell^2(X))$.

\smallskip

\noindent
By replacing $S$ by an isometry with infinite codimension, we obtain a nonunitary isometry in
$\cQ( \calH)$.  It is unclear whether $\cB( \calH)$ or $\cQ( \calH)$ is properly infinite, or whether $\cQ( \calH)$
is simple: it could be that every bounded operator from $\cY$ to $\cY^\perp$ is compact, and if true then
$\cQ( \calH)\cong\cQ(\cY)\oplus\cQ(\cY^\perp)$.
See however~Corollary~\ref{C.XCFstablyfinite} and \S\ref{S.stablyfinite}.

\end{Example}

\subsection{$\cB(\ell^2(X))$ for CF Sets}

Now suppose $X$ is CF.  Then every operator on $\ell^2(X)$ has a lot of finite-dimensional
invariant subspaces:

\paragraph{}
\begin{Lemma}\label{InvSubspaceProp}
Let $X$ be a CF set, $m\geq 1$, and $T_1,\dots,T_m\in\cB(\ell^2(X))$.  Then every finite subset $E$ of $X$ is contained
in a finite subset $F$ such that $\Span\{\xi_x:x\in F\}$ is invariant under $T_1,\dots,T_m$ (by abuse of terminology,
we say $F$ is invariant under $T_1,\dots,T_m$).
\end{Lemma}

\begin{Proof}
Let $E$ be a nonempty finite subset of $X$ (e.g.\ a singleton).  Set $F_0=E$, and inductively
let $F_{n+1}$ be the union of $F_n$ and the supports of $T_k\xi_x$ for all $x\in F_n$ and all $k$ (these supports
are finite by \ref{StrongDFFiniteSupport}).  The $F_n$ form an increasing sequence of finite
subsets of $X$, and hence must stabilize since $X$ is CF.  If $F$ is the final $F_n$,
then $F$ is invariant under $T_1,\dots,T_k$.
\end{Proof}

%\paragraph{}\label{FinitePartitionSAOp}
%Note, however, we cannot conclude (in the absence of Choice) that $X$ decomposes as the disjoint union
%of finite sets invariant under $T$.

\paragraph{}
\begin{Proposition}\label{FinInvSubspaceProp}
Let $X$ be a CF set, and $T_1,\dots,T_m\in\cB(\ell^2(X))$.  Then $X$ partitions into
the disjoint union of finite subsets, each invariant under $T_1,\dots,T_m$.
\end{Proposition}

\begin{Proof}
Replacing each $T_j$ by its real and imaginary parts, we may assume the $T_j$ are self-adjoint.
Define a relation $\sim$ on $X$ by $x\sim y$ if there are $n$
and $x=x_1,x_2,\dots,x_n=y$ such that for each $k$, $1\leq k\leq n-1$, there is a $j$ with
$\langle T_j\xi_{x_k},\xi_{x_{k+1}}\rangle\neq0$ ($j$ can vary with $k$).
This relation is symmetric since $T_j=T_j^*$ for all $j$, transitive, and reflexive (if $T_j\neq0$
for some $j$).  Thus $X$
partitions into $\sim$-equivalence classes.  If $X$ is CF, each equivalence class
is finite by \ref{InvSubspaceProp}.  Each equivalence class is invariant under $T_1,\dots,T_m$.
\end{Proof}

\paragraph{}
\begin{Corollary}\label{C.XCFstablyfinite}
Let $X$ be a CF set.  Then $\cB(\ell^2(X))$ is stably finite.
\end{Corollary}

\begin{Proof}
Let $T$ be an isometry in $\cB(\ell^2(X))$.  Suppose $T$ is not unitary, i.e.\ not surjective, and
let $\eta$ be a vector not in the range of $T$.  Let $E$ be the support of $\eta$, and let $F$ be a
finite set containing $E$ invariant under $T$.  Then the restriction of $T$ to $\cY=\Span\{\xi_x:x\in F\}$
is an isometry in $\cB(\cY)$.  But $\cY$ is finite-dimensional, so $T|_{\cY}$ is unitary,
i.e.\ surjective, contradicting that $\eta\in\cY$ is not in the range of $T$.  Thus $T$ is unitary.
We conclude that $\cB(\ell^2(X))$ is finite.  Since $X\times\{0,\dots,n-1\}$ is also CF
for any $n$, $\cB(\ell^2(X))$ is stably finite.
\end{Proof}

\smallskip

Although a stably finite \cstar-algebra can have an infinite quotient (e.g., the cone over the Toeplitz algebra), we can even conclude that the Calkin algebra $\cQ(\ell^2(X))$ is stably finite (cf.\ \ref{StablyFinCor}).

\subsection{$\cK(\ell^2(X))$ for CF Sets}

\paragraph{}
\begin{Proposition}\label{CompactFinRankProp}
Let $X$ be a CF set.  Then every compact operator on $\ell^2(X)$ has finite rank.
\end{Proposition}

\begin{Proof}
Let $T\in\cK(\ell^2(X))$.  For $\epsilon>0$ set
$$A_\epsilon=\{x\in X:\|T\xi_x\|\geq\epsilon\}\ .$$
Then $A_\epsilon$ is a finite set: if $F$ is a finite subset of $X$ with $\|T-TP_F\|<\epsilon$,
and $x\notin F$, then
$$\|T\xi_x\|=\|T\xi_x-TP_F\xi_x\|\leq\|T-TP_F\|\|\xi_x\|<\epsilon$$
so $A_\epsilon\subseteq F$.  Then $(A_{1/n})$ is an increasing sequence of finite subsets of $X$,
which must stabilize at some finite $A$ since $X$ is CF.  The range of $T$ is contained
in $\Span\{\xi_x:x\in A\}$.
\end{Proof}

\subsection{How Big Is $\cQ(\ell^2(X))$?}\label{S.Q(H)}

\paragraph{}
One might expect that if $X$ is sufficiently DF, then $\cQ( \calH)$ ($ \calH=\ell^2(X)$) might be ``small.''
We could even potentially have $\cQ( \calH)\cong\bbC$, i.e.\ $\cB( \calH)=\cK( \calH)+\bbC1$.
(Note that if $ \calH$ is any infinite-dimensional Hilbert space, the identity operator on $ \calH$ is
not compact since the closed unit ball of $ \calH$ is not totally bounded, cf.\ \ref{InfDimTotBdd}).
This could potentially happen, though, only if $X$ is rigid:

\paragraph{}\label{PermUnitary}
Any permutation $\pi$ of a set $X$ defines a unitary operator $U_\pi$ on $\ell^2(X)$ by
$$U_\pi\left ( \sum_{x\in X} c_x\xi_x\right ) = \sum_{x\in X} c_x\xi_{\pi(x)}\ .$$
If $\pi$ moves infinitely many elements of $X$, then $U_\pi$ is not a compact perturbation
of a scalar.  Thus permutations of $X$ moving infinitely many elements define nontrivial
unitaries in $\cQ(\ell^2(X))$.

\paragraph{}
\begin{Proposition}\label{Q(H)OneDimensional}
	If $X$ is strongly amorphous, then $\cB(\ell^2(X))$ is equal to $\cK(\ell^2(X))+\bbC1$, hence $\cQ(\ell^2(X))$ is one-dimensional and $\cB(\ell^2(X))$ is an AF algebra.
\end{Proposition}

\begin{Proof}
	Let $T\in\cB(\ell^2(X))$.  By \ref{FinInvSubspaceProp}, $X$ partitions into finite subsets
	invariant under $T$.  All but finitely many equivalence classes must be singletons, i.e.\ all but
	finitely many $\xi_x$ are eigenvectors, so $T$ is a compact perturbation of a diagonal operator.
	All but finitely many of the diagonal eigenvalues must be the same since $X$ is amorphous.
	Thus $T$ is a compact (even finite-rank) perturbation of a scalar operator.

	The last sentence follows by Proposition~\ref{P.K(H)simple}.
\end{Proof}

%%%%% Addition

There is an interesting reinterpretation of this result:

\paragraph{}
\begin{Corollary}\label{StrAmHilbSubspaceCor}
	Let $X$ be a strongly amorphous set, and $ \calH=\ell^2(X)$.  Then $ \calH$ is infinite-dimensional,
	and every closed subspace of $ \calH$ has either finite dimension or finite codimension.
\end{Corollary}

\paragraph{}
This is a ``hereditarily indecomposable Hilbert space,'' a Hilbert space version of \cite{GowMa:Unconditional}, where (under AC) a hereditarily
indecomposable Banach space is constructed.  In fact, it is ``better,'' since under AC every
infinite-dimensional Banach space has an infinite-dimensional closed subspace of infinite codimension
(an indecomposable Banach space is just one where no infinite-dimensional closed subspace
has a {\em closed} infinite-dimensional complement).  Our $ \calH$ also has the property that every
compact operator on $ \calH$ has finite rank; under CC every infinite-dimensional Banach space has
an infinite-rank compact operator. (However, by  \cite{argyros2011hereditarily} there is a separable Banach space such that every operator on it is of the form scalar + compact.) Of course, $ \calH$ only exists under set-theoretic assumptions
analysts might find objectionable.

It is possible that ``closed'' can be removed from the statement of \ref{StrAmHilbSubspaceCor}
(note that $ \calH$ has nonclosed subspaces by \ref{CauchyCompleteEx}).  This is just a vector space problem, since $ \calH$
is a complex vector space with a Hamel basis indexed by~$X$.

%%%%% End Addition

\paragraph{} \label{APholds} In the situation of Proposition~\ref{Q(H)OneDimensional}, $\cB(\ell^2(X))$ is also peculiar.
By a classical result of Szankowski (\cite{szankowski1981}), $\cB(\ell^2(\bbN))$ does not have the Approximation Property for Banach spaces. (Recall that a Banach space $X$ has the \emph{Approximation Property} if the identity operator on $X$ can be approximated by finite rank operators uniformly on compact sets.) This implies that if $ \calH$ has a closed subspace isomorphic to $\ell^2(\bbN)$, then $\cB( \calH)$ does not have the Approximation Property. However, if $X$ is strongly amorphous then Propositon~\ref{Q(H)OneDimensional} implies that $\cB(\ell^2(X))$ is AF. In this case conditional expectations to finite-dimensional \cstar-subalgebras witness the
Approximation Property of $\cB(\ell^2(X))$.

\paragraph{}
If $X$ is a strongly amorphous set, and
$Y=X\times\{0,\dots,n-1\}$, then $\cQ(\ell^2(Y))\cong M_n(\bbC)$.  The most interesting case
may be if $X$ is CF but far from rigid, i.e.\ has many nontrivial permutations; then
$\cQ(\ell^2(X))$ is infinite-dimensional and stably finite.  It is hard to see how such a $\cQ$
could have any natural traces.  (Our original motivation for this work was to try to find an example
of a stably finite unital C*-algebra without trace, although it now does not seem promising that such
an example can be constructed by these methods unless we give up the Hahn-Banach Theorem,
without which functional analysis is essentially impossible.)

\smallskip

However, there is no reason to believe that $\cQ( \calH)$ is necessarily simple (in the case where $X$ is DF
of some kind).  Even in the presence of AC, the Calkin algebra associated to a non-separable Hilbert space is not simple: the lattice of closed ideals of $\cQ(\ell^2(\kappa))$ is isomorphic to the poset of infinite cardinals $\leq \kappa$ (see e.g., \cite[\S 12.3.1]{farah2019combinatorial}).  Something similar, but more complicated, may be true for $\cQ(\ell^2(\kappa))$ for $\kappa$ DF (if $\kappa$ is DF,
$\ell^2(\kappa)$ is non-separable).  In the absence of CC, there may even be Hilbert spaces which are so
different that every operator between them is compact. For such pair of spaces,   the Calkin algebra of their direct sum is the direct sum of their Calkin algebras, and therefore non-simple  and even with nontrivial center.  Therefore direct sums of such Hilbert spaces
would have nonsimple Calkin algebras ($\ell^2(X)$, when $X$ is a Russell set, may have this property;  see \ref{RussellOrthSet}).

\smallskip

Actually, this happens.
The following uses a property of Cohen's original CF set (\cite{CohenSet}) that sits between being rigid and strongly rigid and whose proof is included in Proposition~\ref{P.CF.rigid}.

\paragraph{}  \begin{Proposition}\label{P.Calkin.abelian} If $X$ is the CF set that has the property that  for every partition of  $X$  into nonempty finite sets, all but finitely many cells are singletons,  then the Calkin algebra $\cQ(\ell^2(X))$ is a non-separable abelian \cstar-algebra.
\end{Proposition}

\begin{Proof} We first prove that $\cQ(\ell^2(X))$ is equal to its commutative subalgebra,  $\ell^\infty(X)/c_0(X)$.
	Fix $T\in\cB(\ell^2(X))$.  Then Proposition~\ref{FinInvSubspaceProp} implies that~$X$ can be partitioned into finite sets $(X_i)$ such that $\ell^2(X_i)$ is a reducing subspace for $T$ for all $i$.
	By the assumption, all but finitely many of the $X_i$ are
	singletons.  This implies that~$T$  is a compact perturbation of a diagonal operator.  Since $T$ was arbitrary,  $\cQ(\ell^2(X))$ is included in $\ell^\infty(X)/c_0(X)$, and therefore equal to it.

	In order to prove that this algebra is non-separable, we recall that $X$ is a set of reals and consider it with the subspace topology. We first note that $X$ has at most finitely many isolated points,\footnote{As a matter of fact, by construction $X$ has no isolated points.} and we can therefore assume that~$X$ has no isolated points. We may also assume $X$ is bounded, and therefore the closure $\overline X$ of $X$ is a perfect subset of $\R$.
	For every function $f\colon X\to \bbC$, the operator with eigenvalue $f(x)$ corresponding to the eigenvector $\delta_x$ for $x\in X$ belongs to $\ell^\infty(X)/c_0(X)$.
	For $c\in \overline X$  let $f_c$ be the characteristic function of the interval $(-\infty,c)$. If $c<d$, then $X\cap (c,d)$ is infinite, hence  (with $\pi\colon \cB(\ell^2(X))\to \cQ(\ell^2(X))$ we have that $\|\pi(f_c-f_d)\|=1$. Therefore $\{\pi(f_c)\mid c\in (a,b)\}$ is an uncountable discrete subset of $\ell^\infty(X)/c_0(X)$, as required.
\end{Proof}

The following complements Proposition~\ref{P.Calkin.abelian}.

\paragraph{} \begin{Proposition}\label{S.Q.nonabelian}
	If there exists a CF set, then there exists a set $X$ such that $\cB(\ell^2(X))$ is stably finite and $\cQ(\ell^2(X))$ is noncommutative.
\end{Proposition}

\begin{Proof}
	If $Y$ is CF, then so is $X=Y\times Y$ (Proposition~\ref{PFProdProp}), and $\cB(\ell^2(X))$ is stably finite (Corollary~\ref{C.XCFstablyfinite}). We will prove that $\cQ(\ell^2(X))$ contains an isomorphic copy of $M_n(\bbC)$ for every $n\in \bbN$. Fix $F\subseteq Y$ of cardinality~$n$.
	Every permutation $\pi$ of $F$ defines a permutation $\tilde \pi$ of $X$  that agrees with $\pi\times \id_Y$ on $F\times Y$ and fixes all points in $(Y\setminus F)\times X$. This defines a unitary $U_\pi$ on $\cB(\ell^2(X))$. Similarly, every $g\in \ell^\infty(F)$ defines an element of $\ell^\infty(X)$ by $T_g((a,b))=g(a)$ is $a\in F$ and $T_g((a,b))=0$ otherwise. The algebra generated by all $U_\pi$ and all $T_g$ is isomorphic to the subalgebra of $M_n(\bbC)$ generated by all permutation matrices and all diagonal matrices, which is $M_n(\bbC)$. Since $Y$ is infinite, the quotient map acts on this algebra as an isometry, and this concludes the proof.
\end{Proof}

\subsection{Stably Finite \cstar-algebras With no Quasitraces} \label{S.stablyfinite}
The original impetus for this work came from one of the most notorious open problems on \cstar-algebras (see  \cite{haagerup2014quasitraces}): Does every stably finite \cstar-algebra have a tracial state?
In the absence of the Axiom of Choice we provide two examples of stably finite \cstar-algebras without tracial states (and even without quasitraces) of varying levels of satisfactoriness.

\paragraph{ }\label{PropertyOfBaire} Known results provide an easy example of a \cstar-algebra which is stably finite but has no quasitraces in every model from the most studied family of models of ZF in which the Axiom of Choice fails. In the original Solovay's model in which all sets of reals are Lebesgue measurable (\cite{So:Model}), and in other models of the Axiom of Determinacy (see e.g.\cite{kanamori2008higher}), all sets of reals have the Property of Baire. (Notably, unlike the statement that all sets of reals are Lebesgue measurable, the assertion that all sets of reals, and all subsets of any Polish space, have the Property of Baire does not require an inaccessible cardinal. By \cite[\S 7]{shelah1984can}, the consistency of ZFC implies the consistency of ZF+`all  sets of reals have the property of Baire'.)

\paragraph{ }\begin{Lemma} \label{L.Baire} The \cstar-algebra $\ell^\infty(X)/c_0(X)$ has a state if and only if there is a finitely additive probability measure on $\cP(X)$ that vanishes on singletons.
\end{Lemma}

\begin{Proof} Identify $\cP(X)$ with the lattice of projections of $\ell^\infty(X)$.    We claim that finitely additive probability measures on $\cP(X)$ are in a bijective correspondence with states on $\ell^\infty(X)$.  The restriction of a state of the latter  gives a finitely additive probability measure. Conversely, a finitely additive probability measure on $\cP(X)$ can be linearly extended to a state on the $*$-algebra of all finite linear combinations of projections in $\ell^\infty(X)$. This state has norm 1, and  the set of  linear combinations of projections is dense in $\ell^\infty(X)$, and therefore its continuous extension is a state on $\ell^\infty(X)$.

	Finally, a finitely additive measure on $\cP(X)$ vanishes on the singletons if and only if the associated state vanishes on $c_0(X)$, and this completes the proof.
\end{Proof}

\paragraph{ }\begin{Corollary}\label{C.Baire} If all sets of reals have the Property of Baire then the commutative \cstar-algebra $\ell^\infty(\N)/c_0(\N)$ has no nontrivial bounded linear functionals.

	In particular, it does not have a nontrivial representation on a Hilbert space and it is not isomorphic to $C(X)$ for a compact Hausdorff space $X$.
\end{Corollary}

\begin{Proof}
	As pointed out without a proof in \cite[p. 3]{So:Model}, if all sets of reals have the Property of Baire then there is no finitely additive probability measure on $\cP(\N)$ that vanishes on singletons (for a proof see e.g., \cite[(8) on p. 206]{pincus2006strength}).
	This implies that the commutative \cstar-algebra $A=\ell^\infty(\N)/c_0(\N)$ has no nontrivial bounded linear functionals.

	Suppose that $\pi \colon A\to\cB( \calH)$ is a nontrivial representation on some Hilbert space $ \calH$. Then there are $\xi\in  \calH$ and $T$ in the image such that $\langle T\xi,\xi\rangle\neq 0$, and therefore $\varphi(a)=\langle \pi(a)\xi,\xi\rangle$ is a nontrivial state on $A$; contradiction.

	Finally, suppose that $A$ is isomorphic to $C(X)$ for a compact Hausdorff space $X$, or even that $\Phi\colon A\to C(X)$ is a nontrivial *-homomorphism. Then for some $x\in X$ the composition of $\Phi$ with the evaluation functional on $C(X)$ is a nontrivial state on $A$; contradiction.
\end{Proof}

Corollary~\ref{C.Baire} is not satisfactory, since the reason for the absence of tracial states in $\ell^\infty(\N)/c_0(\N)$ is the failure of the Hahn--Banach theorem.
The assumption that all sets of reals have the Property of Baire implies that $\cB(\ell^2(\N))$ has a unique pure state, up to the unitary equivalence (\cite[Exercise~12.6.6]{farah2019combinatorial}), making it a counterexample to Naimark’s Problem (see \cite{AkeWe:Consistency} and \cite[\S 11.2]{farah2019combinatorial}).

We don't know whether the assumption of the following proposition is relatively consistent with ZF (see \ref{P.CF.no-measure}).

\paragraph{} \begin{Proposition} \label{P.stably-finite-no-tracial-states} Suppose that there is a CF set~$X$ such that $\cP(X)$  has no finitely additive probability measure that vanishes on singletons.
	Then  $\cB(\ell^2(X))$ is stably finite, it has no tracial states, but its states separate the points.
\end{Proposition}

\begin{Proof}   Since $X$ is CF, $\cB(\ell^2(X))$ is stably finite by  Corollary~\ref{C.XCFstablyfinite}.
	Suppose~$\tau$ is a tracial state of $\cB(\ell^2(X))$. We claim that $\tau$ vanishes on compact operators. Otherwise, since projections of finite rank form an approximate unit in $\cK(\ell^2(X))$, we can fix a projection~$p$ of rank 1 and $n\geq 1$ such that $\tau(p)>1/n$. Since $X$ is infinite, we can find $n$ orthogonal projections of rank 1 in $\cB(\ell^2(X))$, $p_1,\dots, p_n$. All of these projections are Murray--von Neumann equivalent, and therefore $\tau(p_j)=\tau(p)>1/n$ for all $j$. This implies $\tau(\sum_j p_j)>1$, contradicting the assumption that $\tau$ is a state.

	Therefore the restriction of $\tau$ to $\ell^\infty(X)$ defines a state on $\ell^\infty(X)/c_0(X)$. Since no finitely additive measure on $\cP(X)$ vanishes on singletons, we have contradiction by Lemma~\ref{L.Baire}.

	Finally, any two distinct operators in $\cB(\ell^2(X))$ are separated by a vector state associated with some  $\xi\in \ell^2(X)$, and this completes the proof.
\end{Proof}

\subsection{Spectrum of Operators}\label{S.Spectrum}

\paragraph{}
If $X$ is sufficiently DF, one might expect operators on $ \calH=\ell^2(X)$ to have small spectrum.
For example, if $X$ is CF, \ref{CompactFinRankProp} implies that the spectrum of every compact operator on
$ \calH$ is finite.

\paragraph{}
\begin{Theorem}\label{T.PDF-spectrum}
Let $X$ be a power Dedekind-finite set, and $T\in\cB(\ell^2(X))$, Then $\sigma(T)$ is finite,
and each $\lambda\in\sigma(T)$ is an eigenvalue.
\end{Theorem}

\begin{Proof}
We may assume $X$ is infinite (PF) since the finite case  is well known.  Thus $X$ is CF,
so by \ref{FinInvSubspaceProp} $X$ partitions into finite subsets $\{X_i:i\in I\}$ whose spans $ \calH_i$ are invariant under $T$.
There is a finite-to-one function $f:X\to I$.  Let $Y\subseteq\bbC$ be $\cup_{i\in I}S_i$,
where $S_i$ is the (finite) set of eigenvalues of the restriction of $T$ to $ \calH_i$.

\smallskip

\noindent
We claim $Y$ is finite.  If not, then there is a partition of $Y$ into a disjoint sequence of
nonempty subsets $Y_n$ (\ref{DedInfRSubProp}).  For each $n$ let $I_n$ be the set of all $i\in I$ such that
$S_i\cap Y_n\neq\emptyset$.  Since $S_i$ is finite, each $i\in I$ is in only finitely many $I_n$.  But $\cup_nI_n=I$,
an infinite set, so there are an infinite number of distinct $I_n$.  If $A_n=f^{-1}(I_n)$, then
each $A_n$ is a subset of $X$, and infinitely many are distinct, contradicting that
$X$ is PF.

\smallskip

\noindent
Now we claim that $Y$ is the entire spectrum of $T$.  Let $\lambda\in\bbC\setminus Y$.
Then $T-\lambda 1$ is bounded below on $ \calH_i$ for each $i\in I$.  For each $n$, let $J_n$
be the set of $i\in I$ such that $T-\lambda 1$ is bounded below by $\frac{1}{n}$ on $ \calH_i$,
and let $B_n=f^{-1}(J_n)\subseteq X$.  We have $J_n\subseteq J_{n+1}$ and hence
$B_n\subseteq B_{n+1}$ for all $n$; $\cup_nJ_n=I$ so $\cup_n B_n=X$.  Since $X$ is PF,
we must have $B_n=X$ for some $n$, i.e.\ $T-\lambda 1$ is bounded below by $\frac{1}{n}$
on all of $\ell^2(X)$ and hence invertible.
\end{Proof}

The conclusion of Theorem~\ref{T.PDF-spectrum} does {\em not} imply that every operator is the sum of a scalar and a compact operator. For example,
if  $X=\bigsqcup_{z\in Z} X_z$, where $Z$ is strongly amorphous, is  a set as in \S\ref{SuperRussell}, then the permutation that swaps each of the pairs $X_z$ moves all points of $X$, and therefore defines a unitary which is not a compact perturbation of a scalar (\S\ref{PermUnitary}).

\paragraph{}\label{LargeSpectrum}
If $X$ is weakly Dedekind-infinite, there are operators on $ \calH=\ell^2(X)$ with large spectrum.  Write $X$
as a countable union of pairwise disjoint nonempty sets $A_n$.  If $K$ is a separable compact
subset of $\bbC$ (separability of subsets of $\R$ or $\bbC$ is not automatic without CC,
but there are many separable compact subsets, e.g.\ the closure of any countable bounded
subset), then there is a diagonal operator $T$ on $ \calH$ whose spectrum is $K$, defined by taking a dense sequence
$(a_n)$ in $K$, for $x\in X$ setting $n(x)$ the $n$ for which $x\in A_n$, and setting
$$T\left ( \sum_{x\in X} c_x\xi_x\right ) =\sum_{x\in X}a_{n(x)}c_x\xi_x\ .$$
If all the $A_n$ are infinite (which can be arranged by combining them), the spectrum of the image
of $T$ in $\cQ( \calH)$ is also $K$.

\paragraph{}
Let $X$ be a CF set which is not PF, and write $X$ as a disjoint union of a sequence $(A_n)$
of nonempty subsets.  Let $T$ be the operator on $\ell^2(X)$ which is $\frac{1}{n}$ times the
identity on the span of $A_n$.  Then $T$ is bounded ($\|T\|=1$) and injective, and is also
surjective since every vector in $\ell^2(X)$ has finite support (Proposition~\ref{StrongDFFiniteSupport}).  But the inverse $T^{-1}$ is
not bounded.  Thus the Bounded Inverse Theorem, Open Mapping Theorem, and Closed Graph
Theorem all fail in this situation.\footnote{It is known that the Open Mapping Theorem is roughly equivalent to the countable Axiom of Choice, \cite{fellhauer2017relation}.}  By considering the finite truncations of $T^{-1}$, a sequence of
pointwise bounded but not uniformly bounded operators is obtained, showing that the Uniform
Boundedness Theorem also fails.  The usual proofs of these theorems use some version of the Baire Category
Theorem, equivalent to the Axiom of Dependent Choice, which implies CC and thus is incompatible with
the existence of DF sets.

There is another similar example:  the Banach spaces $\ell^p(X)$, $1\leq p<\infty$, are all the
same set of functions, and the norms are comparable but pairwise nonequivalent.

\paragraph{}
It would also be interesting to study the properties of $\cB( \calH)$, etc., for Hilbert spaces $ \calH$
which do not have orthonormal bases.  It is difficult to find tools to analyze these cases, however.
Even a Hilbert space $ \calH$ with an orthonormal basis may have closed subspaces without one (\ref{L.ExistenceOfBasis} (2)),
hence $\cB( \calH)$ can have corners without sufficiently many rank 1 projections.
%One question which may help: is every Hilbert space isometric to a closed subspace of a Hilbert
%space with an orthonormal basis?

\section{New Size and Comparison Relations on Hilbert Spaces and Sets}\label{S.HDF}

%HDF, HCF, HPF, Hilbert-amorphous; HDF=HCF=CF; finite and stably finite $\cB(H)$ and $\cQ( \calH)$

\subsection{Hilbert Dedekind-finiteness}

There are Hilbert space analogs of the principal flavors of Dedekind-finiteness, which lead to
potentially new flavors on the set-theoretic level:

\paragraph{}
\begin{Definition}\label{HilbDFDef}
Let $ \calH$ be a Hilbert space.  Then
\begin{enumerate}
\item[(i)]  $ \calH$ is {\em Hilbert Dedekind-infinite} if $ \calH$ is isometrically isomorphic to a proper
subspace of itself.  Otherwise $ \calH$ is {\em Hilbert Dedekind-finite}.
\item[(ii)]  $ \calH$ is {\em Hilbert Cohen-infinite} if $ \calH$ has a sequence of mutually orthogonal
nonzero finite-dimensional subspaces.  Otherwise $ \calH$ is {\em Hilbert Cohen-finite}.
\item[(iii)]  $ \calH$ is {\em Hilbert dually Dedekind-infinite} if $\cB( \calH)$ contains an operator
which is surjective but not injective.  Otherwise $ \calH$ is {\em Hilbert dually Dedekind-finite}.
\item[(iv)]  $ \calH$ is {\em Hilbert weakly Dedekind-infinite} if $ \calH$ has a sequence of mutually orthogonal
nonzero subspaces.  Otherwise $ \calH$ is {\em Hilbert Power Dedekind-finite}.
\item[(v)]  $ \calH$ is {\em Hilbert-amorphous} if $ \calH$ is infinite-dimensional and does not contain two orthogonal
infinite-dimensional subspaces, i.e.\ every closed subspace either has finite dimension or
finite codimension.
\end{enumerate}
If $X$ is a set, say $X$ is Hilbert Dedekind-finite, etc., if $\ell^2(X)$ is Hilbert Dedekind-finite, etc.
\end{Definition}

\subsubsection*{Notation} For a Hilbert space, HDF (resp.\ HCF, HDDF, HPF) means ``in\-fi\-ni\-te-dimensional and Hilbert Dedekind-finite
(resp. Hilbert Cohen-finite, Hilbert dually Dedekind-finite, Hilbert power Dedekind-finite),
and similarly for a set.

\paragraph{}
It is clear that Hilbert-amorphous $\imp$ HPF $\imp$ HCF $\imp$ HDF (cf.\ \ref{HDFThm}), and that
for sets HDF $\imp$ DF, HCF $\imp$ CF, HPF $\imp$ PF, and Hilbert-amorphous $\imp$
amorphous (the usual notions are the special cases of the Hilbert versions for subspaces
spanned by sets of standard basis vectors).
Also, $ \calH$ is Hilbert Dedekind-finite if and only if every isometry in $\cB( \calH)$ is unitary, i.e.\ if and only if $\cB( \calH)$ is finite.

\bigskip

There is considerable collapsing among these notions, particularly the set ones:

\paragraph{}
\begin{Proposition}\label{HDFThm}
Let $ \calH$ be a Hilbert space.  The following are equivalent:
\begin{enumerate}
\item[(i)]  $ \calH$ contains an infinite orthonormal sequence.
\item[(ii)]  $ \calH$ is Hilbert Dedekind-infinite.
\item[(iii)]  $ \calH$ is Hilbert dually Dedekind-infinite.
\end{enumerate}
\end{Proposition}
\begin{Proof}
(i) $\imp$ (ii):  If $ \calH$ contains an orthonormal sequence, then $ \calH$ contains a closed subspace isometrically
isomorphic to $\ell^2(\N)$, and there is a nonunitary isometry on this subspace, which can be
expanded to a nonunitary isometry on $ \calH$ by making it the identity on the orthogonal complement.

\smallskip

\noindent
(ii) $\imp$ (i):  if $V$ is a nonunitary isometry in $\cB( \calH)$, and $\cY$ is the range of $V$,
let $\xi$ be a unit vector in $\cY^\perp$.  Then $\xi$, $V\xi$, $V^2\xi$, $\dots$ is an orthonormal sequence in $ \calH$.

\smallskip

\noindent
(ii) $\imp$ (iii):  If $V$ is a nonunitary isometry in $\cB( \calH)$, then $V^*$ is surjective but
not injective.

\smallskip

\noindent
(iii) $\imp$ (ii):  If $T$ is surjective but not injective, let $\cY=\cN(T)^\perp$.  Then the restriction
of $T$ to $\cY$ is an injective map from $\cY$ to $ \calH$, so if $T=V|T|$ is the polar decomposition of $T$, then $V$ is a coisometry,
and $V^*$ is an isometry from $ \calH$ onto $\cY$.
\end{Proof}

\paragraph{}
\begin{Corollary}\label{HDFCor}
Let $X$ be an infinite set.  The following are equivalent:
\begin{enumerate}
\item[(i)]  $X$ is CF.
\item[(ii)]  $X$ is HDF.
\item[(iii)]  $X$ is HDDF.
%\item[(iv)]  $X$ is HCF.
\end{enumerate}
\end{Corollary}
\begin{Proof}
By \ref{HDFThm}, (ii) $\ifff$ (iii).

(ii) $\ifff$ (i): By \ref{HDFThm}, $X$ is not HDF if and only if  $\ell^2(X)$ contains a countably infinite orthonormal sequence. By \ref{L.SDFBasis}, this is equivalent to $X$ not being CF.   \end{Proof}

\paragraph{}\label{HCFandsuch}
It is unclear whether HCF is equivalent to CF (and hence to HDF and HDDF) for sets (but see \ref{StrongDFFiniteSupport} and \ref{L.SDFBasis}).
HPF is distinct from the conditions of \ref{HDFCor} for sets since HPF $\imp$
PF and there are sets which are CF but not PF.  It is unclear whether HPF is the same as PF
since closed subspaces of $\ell^2(X)$, $X$ PF, are potentially skew from the standard basis.
Thus HPF is potentially a Dedekind-finiteness condition on sets strictly more restrictive than PF,
and HCF is also a potentially new condition.

\paragraph{}
Similarly, Hilbert-amorphous is potentially strictly stronger than amorphous for sets.  However, by Corollary~\ref{StrAmHilbSubspaceCor} a strongly amorphous set is Hilbert-amorphous.  The converse is quite unclear.  Thus
Hilbert-amorphous is potentially a Dedekind-finite condition strictly between being amorphous
and strictly amorphous.

Here is an extended version of Fig.~\ref{Fig.DF} from Section~\ref{S.RussellSets}.

\bigskip

A simple but interesting observation is:

\paragraph{}
\begin{Proposition}\label{HDFSumProp}
A finite direct sum of HDF Hilbert spaces is HDF.
\end{Proposition}

\begin{Proof}
Let $ \calH_1,\dots, \calH_n$ be Hilbert spaces.  Suppose $(\xi_j)$ is an orthonormal sequence
in $ \calH= \calH_1\oplus\cdots\oplus \calH_n$.  Then the subspace of $ \calH$ spanned by the $\xi_j$
is infinite-dimensional.  For each $k\leq n$, let $(\xi_{j,k})$ be the sequence of $k$'th coordinates
of $\xi_j$, and let $\cX_k$ the closed subspace of $ \calH_k$ spanned by the $\xi_{j,k}$.
The $\cX_k$ cannot all be finite-dimensional, so there is a $k$ such that the sequence $(\xi_{j,k})$
contains a linearly independent subsequence, from which an orthonormal sequence can be made by Gram--Schmidt.  Thus $ \calH_k$ is not HDF.
\end{Proof}

\paragraph{}
\begin{Proposition}\label{StablyFinCor}\label{QFiniteCor}

	Let $ \calH$ be a Hilbert space.  Then the following are equivalent:
	\begin{enumerate}
		\item[(i)]  $\cB( \calH)$ is finite.
		\item[(ii)]  $\cB( \calH)$ is stably finite.
		\item[(iii)]  $\cQ( \calH)$ is finite.
		\item[(iv)]  $\cQ( \calH)$ is stably finite.
	\end{enumerate}
\end{Proposition}

\begin{Proof}	Since $M_n(\cB(H))$ is isomorphic to $\cB(H^n)$, (i) $\ifff$ (ii) follows from \ref{HDFSumProp}.

(i) $\ifff$ (iii):	If $\cB( \calH)$ is infinite, then by \ref{HDFThm} $ \calH$ contains an orthonormal sequence, hence an isometry of
	infinite codimension, whose image in $\cQ( \calH)$ is a nonunitary isometry, so $\cQ( \calH)$ is
	also infinite.  Conversely, if $v$ is a nonunitary isometry in $\cQ( \calH)$, let $T\in\cB( \calH)$
	have image $v$.  If $T=V|T|$ is the polar decomposition, then the image of $V$ is also $v$.
	$V$ is a partial isometry whose source projection $V^*V$ has finite codimension and range
	projection $VV^*$ has infinite codimension.  There is thus a partial isometry $W\in\cB( \calH)$
	with $W^*W=1-V^*V$ and $WW^*\leq 1-VV^*$ (any infinite-dimensional Hilbert space
	has closed subspaces of any finite dimension), and $V+W$ is a nonunitary isometry in $\cB( \calH)$,
	so $\cB( \calH)$ is also infinite.

This also proves (ii) $\ifff$ (iv) and completes the proof.
\end{Proof}

\subsection{Quantum Cardinals}
Here is another interpretation of our setup.
Consider the categories $\Set$ of all sets with injective maps as morphisms and $\Hilbert$ of all Hilbert spaces with isometries as morphisms. Then  $F(X)= \ell^2(X)$ is a functor between these categories. The AC implies the following:
\begin{enumerate}
	\item [(i)] $F$ is a bijection on the objects (although it is not an equivalence of categories, since Hilbert spaces have many more morphisms).
	\item [(ii)] $|X|\leq |Y|$ if and only if $\ell^2(X)$ is isomorphic to a closed subspace of~$\ell^2(Y)$.
\end{enumerate}
One of the facets of our project is the study of this functor in the context when the AC fails.
Naturally separating  (i) into the statements `$F$ is surjective' and `$F$ is injective', we will see that
both can fail: There can be  a Hilbert space without a basis (Theorem~\ref{L.ExistenceOfBasis})
and a Hilbert space can have orthonormal bases of different cardinalities (Example~\ref{NonuniqueBasis}).

In order to discuss (ii), we consider the following, potentially new, relation between cardinals based on Hilbert spaces.

\paragraph{}
\begin{Definition}\label{PreceqDef}
	Let $X$ and $Y$ be sets, and $\kappa=|X|$, $\lambda=|Y|$.  We write
	\[
	\lambda\preceq\kappa
	\]
	if $\ell^2(X)$ contains a (closed) subspace isometrically isomorphic to $\ell^2(Y)$, or equivalently
	(\ref{OrthBasisL2}) $\ell^2(X)$ contains an orthonormal set of cardinality $\lambda$.
\end{Definition}

\paragraph{}
The relation $\preceq$ is clearly transitive, and $\lambda\leq\kappa\imp\lambda\preceq\kappa$.
But $\preceq$ is much weaker than $\leq$:

\paragraph{}
\begin{Examples}\label{PreceqExs}
	(i)  If $X$ is DF but not CF, then $\aleph_0\preceq|X|$ and $\aleph_0\not\leq |X|$ (Example~\ref{NonuniqueBasis}).

	\smallskip

	\noindent
	(ii)  If $X$ is DF but not CF, then $|X|\pm n\preceq|X|$ for any $n$ (the extra finite-dimensional
	space can be absorbed into an $\ell^2(\N)$ summand, see \S\ref{RussellNotCF}).
	\smallskip

	\noindent
	(iii) By Theorem~\ref{SepHilbSpaceThm}, $|X|\preceq \aleph_0$ if and only if $|X|\leq \aleph_0$. By Theorem~\ref{L.l2(kappa)}, if $\kappa$ is a well-ordered cardinal then $|X|\preceq\kappa$ if and only if $|X|\leq \kappa$.
\end{Examples}

The relation $\preceq$ is not antisymmetric (\ref{PreceqExs}(ii)), but there is a partial result.  We have the following Hilbert
space version of the Schr\"{o}der-Bernstein Theorem; the proof (in ZF) is a simple modification of
Bernstein's proof of the set-theoretic result, and this Hilbert space result may have been
previously noted (although under AC it is a trivial consequence of the existence and uniqueness
of orthogonal dimension).

\paragraph{}
\begin{Theorem}\label{HilbSchBernThm}
	Let $ \calH$ and $ \calH'$ be Hilbert spaces.  If each is isometrically isomorphic to a (closed) subspace
	of the other, then $ \calH$ and $ \calH'$ are isometrically isomorphic.

In particular, for all sets $X$ and $Y$, if $|X|\preceq|Y|$ and $|Y|\preceq|X|$, then $\ell^2(X)\cong\ell^2(Y)$.
\end{Theorem}

\begin{Proof}
	Let $T: \calH\to \calH'$ and $S: \calH'\to \calH$ be isometries (not necessarily surjective).  Let
	$\cX_1=[S( \calH')]^\perp$ and $\cY_1=T(\cX_1)$, and recursively define $\cX_{n+1}=S(\cY_n)$
	and $\cY_{n+1}=T(\cX_{n+1})$ for each $n$.  The $\cX_n$ are mutually orthogonal subspaces of $ \calH$, and the $\cY_n$ mutually orthogonal subspaces of~$ \calH'$.

	\smallskip

	\noindent
	Let $\cX$ be the closed span of $\cup_n\cX_n$
	and $\cY$ the closed span of $\cup_n\cY_n$.  Then~$T$ maps $\cX$ isometrically onto $\cY$.
	We claim that $S$ maps $\cY^\perp$ isometrically onto~$\cX^\perp$.  If $\cZ$ is the closed span
	of $\cup_{n=2}^\infty\cX_n$, then $S$ maps $\cY$ isometrically onto $\cZ$, and so maps
	$\cY^\perp$ onto $\cX_1^\perp\cap\cZ^\perp=\cX^\perp$.

	\smallskip

	\noindent
	Let $U: \calH\to \calH'$ be equal to $T$ on $\cX$ and $S^{-1}$ on $\cX^\perp$.  Then $U$ is an isometry
	from $ \calH$ onto $ \calH'$.
\end{Proof}

\paragraph{}
Thus we also obtain a new equivalence relation $\sim$ on cardinals, where $|X|\sim|Y|$ if
$\ell^2(X)\cong\ell^2(Y)$.

\paragraph{}
For a partial characterization, say a set $X$ {\em finitely covers} a set $Y$ if there is a
finite-to-one function from $X$ onto $Y$.  By the argument of
\ref{StrongDFFiniteSupport}[(iii) $\imp$ (i)], if a subset of a set $X$ finitely covers a set $Y$,
then $|Y|\preceq|X|$.  The converse, however, seems doubtful, and we do not know a purely
set-theoretic characterization of either $\preceq$ or $\sim$.  It may be that for sufficiently
Dedekind-finite cardinals, $\preceq$ coincides with $\leq$, and thus $\sim$ is just equality.

\section{Models of ZF}\label{S.ModelsofZF}

In this section we collect set-theoretic results used in the earlier sections. The following is taken from \cite{CohenSet}.

\paragraph{} \begin{Proposition}\label{P.Cohen} If ZF is consistent, then so are the theories ZF+`There exists a DF set $X$' and   ZF+`There exists a CF set $X$'.
\end{Proposition}

We  proceed to describe constructions of models of ZF used in previous sections that cannot be found in the literature.
By using the method of \cite{karagila2019iterated}, one can construct a single model of ZF in which all instances of the failure of the Axiom of Choice used in the present paper appear.
First we prove the result used in \S\ref{SuperRussell}, after recalling some definitions. Given a discrete structure $A$, the \emph{definable} subsets are the subsets of $A$ definable by a first-order formula, possibly with parameters in $A$. Thus these are the sets of the form
\[
\varphi^{A,\bar b}=\{a\in A\mid A\models \varphi(a,\bar b)\}
\]
where $\varphi(x,\bar y)$ is a first-order formula in the signature of $A$ and $\bar b$ is a tuple in $A$ of the same sort as $\bar y$. Definable subsets of $A^k$ for $k\geq 2$ are defined analogously. The family of definable subsets of $A$ or $A^k$ for $k\geq 2$ does not depend on the extent of choice available, because each $\varphi^{A,\bar b}$ corresponds to a pair $\varphi,\bar b$.

\paragraph{} \begin{Proposition}\label{L.Truss} If ZF is consistent, then so is the theory ZF+`There exists an amorphous set $X$ which can be presented as $X=\bigsqcup_{z\in Z} X_z$,  where~$Z$ is strongly amorphous, and $X_z$ are pairwise disjoint two-element sets.'
\end{Proposition}

\begin{Proof} We will use the ideas and terminology from \cite[\S 3]{karagila2022choiceless}. There, Theorem~3.2 allows us to take a structure from a model of ZFC and create a model of ZF in which there is a copy of the structure where all the subsets are definable.  So our goal shifts to finding an infinite structure $A$ with the following properties.
	\begin{enumerate}
		\item The definable subsets of $A$ are finite or cofinite,
		\item $A$ can be partitioned into pairs (i.e., there is an equivalence relation on $A$ such that each equivalence class has exactly two elements), and
		\item any partition of the equivalence classes must be almost entirely into singletons.
	\end{enumerate}
In addition, the requirements of \cite[Theorem~3.2]{karagila2022choiceless} include identifying a group of automorphisms of the structure and an ideal of subsets, which in our case will be the ideal of finite subsets.

  Let $R$ be the equivalence relation on $\N$ defined by $R(m,n)$ if and only if $\lfloor m/2\rfloor = \lfloor n/2\rfloor$. In other words, $R$ is the equivalence relation whose classes are the sets $\{2n,2n+1\}$ for $n\in\N$. Let $\cG$ be the group of automorphisms of $(\N,R)$, we will say that a subgroup of $\calH$ is \emph{large} if there is some $m\in\N$ such that the group \[\calH_m=\{\pi\in\cG\mid \pi(n)=n\text{ for all }n<2m\}\] is included in $\calH$. In other words, if the subgroup fixes pointwise the first $m$ equivalence classes. The family of large subgroups is closed under subgroups and intersection, as well as under conjugation (by automorphisms from $\cG$). In the terminology of \cite[\S 2.2]{karagila2022choiceless}, the large groups form a \emph{normal filter of subgroups} of $\cG$.

  For any fixed $k\geq 1$ and $\pi\in\cG$, $\pi^k$ is the naturally defined automorphism acting on $\N^k$. We will say that $A\subseteq\N^k$ is \textit{stable under a subgroup $\calH$ of $\cG$} if $\pi^k[A]=A$ for all $\pi\in\calH$, and we say that $A$ is \emph{stable} if there is some large subgroup under which it is stable.

  The following observations are immediate.
  \begin{enumerate}[resume]
  \item\label{1.ZF}  If $A\subseteq\N$ is stable then it is finite or cofinite.
  \item \label{2.ZF} If in addition $A$ is a union of $R$-equivalence classes, then $A/R$ is finite or cofinite (in $\N/R$).
  \item \label{3.ZF} Moreover, suppose that $E$ is an equivalence relation on $\N$ which is definable over $(\N,R)$ (with parameters, perhaps) such that $E$ is stable and $R\subseteq E$, then there is some $m$ such that for $n\geq 2m$, the $E$-equivalence class of $2n$ is either $\{2n,2n+1\}$ or $\N\setminus\{0,\dots,2m-1\}$. In particular the equivalence relation $E$ induces on $\N/R$ has only finitely many non-singleton equivalence classes.
  \end{enumerate}
  Given a countable transitive model of ZF and applying  \cite[Theorem~3.2]{karagila2022choiceless} we obtain a forcing extension and an intermediate model, $V\subseteq W\subseteq V[G]$ such that in $W$ there is a structure $(X,S)$ for which in $V[G]$ there is an isomorphism $(X,S)\cong(\N,R)$ and any $A\subseteq X$ such that $A\in W$ is the image of some stable subset under the isomorphism. Clearly, $X$ is a union of unordered pairs in $W$. By \eqref{1.ZF}, $X$ must be amorphous in $W$. By \eqref{2.ZF}, $Z=X/S$ is amorphous as well, and by \eqref{3.ZF} it is in fact strongly amorphous, as any partition of $Z$ can be made into a partition of $X$.  Therefore the equivalence relation it defines on $X$ must be one which extends~$S$.
\end{Proof}

Proposition~\ref{P.CF.rigid} below completes the proof of Proposition~\ref{P.Calkin.abelian}.
We briefly recall the basic Cohen’s model (\cite[\S 5.3]{jech2008axiom}) with a CF set. Let $\bbP$ be the forcing for adding side-by-side Cohen subsets of $\bbN$,  $\dot x_n$, for $n\in \bbN$. The conditions in $\bbP$ are the finite partial functions from $\bbN^2$ into $\{0,1\}$, and $\bbP$ is ordered by extension.  If $G\subseteq \bbP$ is sufficiently generic, then $\bigcup G$ is a function from $\bbN^2$ into $\{0,1\}$,   and the interpretation of $m$ is defined by  $\dot x_m^{\dot G}(n)=\dot G(m,n)$ for all $m$ and $n$   (see e.g., \cite[2(a)]{karagila2020have} for additional details). As customary,  $\dot x_m^{\dot G}$ is identified with the subset of $\bbN$ that it is the characteristic function of.

Let $\cG$ denote  the group of all finitary permutations of $\bbN$. Consider the action of~$\cG$ on  $\bbP$ defined as follows. A permutation $g\in \cG$ sends $s\in \bbP$ to $g.s\in \bbP$, where
\[
\dom(g.s)=\{(g(m),n)\mid (m,n)\in \dom(s)\}
\]
and $g.s(g(m),n)=s(m,n)$ for all $(m,n)\in \dom(s)$. This action naturally extends to
$L^\bbP$ (the set of $\bbP$-names in G\"odel's constructible universe $L$; any other model of ZFC would do in place of $L$), so that  $g.\dot x_m=\dot x_{g.m}$  for all $g$ and $m$.

For each $m\in \bbN$, let $\cG_m$ be the subgroup of $\cG$ consisting of permutations~$g$ that fix all $n\leq m$. Then for all $m$ and $g\in \cG$ there exists $n$ such that $g\cG_m g^{-1}\subseteq \cG_n$. Therefore these groups form a filter of subgroups of $\cG$.  Consider the set of all \emph{symmetric} names $L^{\bbP}$-names:
\[
\Sym=\{\dot a\in L^{\bbP}\mid (\exists m) g.\dot a=\dot a\text{ for all }g\in \cG_m\}.
\]
Recursively define the set of \emph{hereditarily symmetric names}, $\HSym$, to  be the set of all names $\dot x$ such that $\dot x$ is forced to be a subset of $\Sym$ and it belongs to $\Sym$. Note that this includes canonical names for all elements of $L$.
If $\dot G\subseteq \bbP$ is generic over $L$, then \cite[Theorem 5.19]{jech2008axiom} implies that the set of all $\dot G$-evaluations of hereditarily symmetric names
\[
M=\{\dot a^{\dot G}\mid \dot a\in \HSym\}
\]
is a model of ZF that contains $X=\{\dot x_m^{\dot G}\mid m\in \bbN \}$, and $X$ is CF in this model.

The following folklore proof was communicated to A.K. by Andreas Blass.

\paragraph{}\begin{Proposition} \label{P.CF.rigid} With the notation from the previous paragraph, the set $X$ has the property that  for every partition of  $X$ in $M$ into nonempty finite sets, all but finitely many cells are singletons.
\end{Proposition}

\begin{Proof}
	Suppose that $\dot\pi$ is a $\bbP$-name, and that some condition $p\in \bbP$ forces that it  is a name for a partition of $X$ into nonempty finite subsets that belongs to $N$. In this proof we write $x_n=\dot x_n^{\dot G}$ and $\pi=\dot\pi^{\dot G}$. 	It will be convenient to write $i\sim j$ if $x_i$ and $x_j$ belong to the same piece of $\pi$.   Fix $m$ such that $g.\pi=\pi$ for all $g\in \cG_m$. Since all pieces of $\pi$ are finite, we can find an extension  $q\leq p$ that decides $k\geq m$ such that for all $i\in\bbN$, if some extension of $q$ forces that $i\sim j$ for some $j<m$, then $i<k$.

	We claim that $q$ forces that for $k\leq \min(i,j)$, $i\sim j$ implies $i=j$.
	Assume otherwise, and fix a condition $r\leq q$ that decides $k<i<j$ such that $i\sim j$. Since the support of $r$ is finite, we can find permutations $g_l\in \cG_m$, for $l\in \bbN$, such that $g_l(i)=i$, $g_l(j)\neq g_{l'}(j)$ when $l\neq l'$, and the supports of the conditions $g_l.r$, for $l\in \bbN$, form a $\Delta$-system, and are pairwise compatible. Note that $r_l$ forces $i\sim g_l(j)$ for all $l$. By genericity, $r$ forces that the set $\{l\mid r_l\in \dot G\}$ is infinite, and therefore the piece of the partition $\pi$ to which $i$ belongs is infinite; contradiction.
\end{Proof}

\section{Problems} \label{S.Problems}
We conclude with a short list of select open problems. Paragraph \ref{HCFandsuch}
 contains several interesting problems that we will not repeat here.

The following question gave the original impetus to this project.

\paragraph{} \begin{Question}\label{Q.stablyfinite.notraces}
Is the existence of  a Hilbert space $ \calH$ such that $\cB( \calH)$ is stably finite but it has no tracial states relatively consistent with ZF?
\end{Question}

By \ref{C.Baire}, the existence of an abstract unital \cstar-algebra that is stably finite but has no tracial states is relatively consistent with ZF. However these algebras have no states (and therefore no nontrivial representations on a Hilbert space) at all. On the other hand, $\cB(\calH)$ is already represented as a concrete \cstar-algebra and a positive answer to Question~\ref{Q.stablyfinite.notraces} would be more interesting.

 A positive answer to the following would imply that the assumption of Proposition~\ref{P.stably-finite-no-tracial-states} is relatively consistent with ZF, and that Question~\ref{Q.stablyfinite.notraces} has a positive answer.

\paragraph{} \begin{Question} \label{P.CF.no-measure} Is the existence of a CF set $X$ such that  no finitely additive probability measure on $\cP(X)$ vanishes on singletons relatively consistent with ZF?
\end{Question}

By the main result of  \cite{pincus1977definability}, there is a model of ZF in which for every infinite set  $X$,  no finitely additive probability measure on $\cP(X)$ vanishes on singletons. However, this model satisfies the Axiom of Dependent Choices (DC) and therefore has no DF sets.

\paragraph{}\begin{Problem}
	Study von Neumann algebras in ZF.
	\end{Problem}

There is an ample supply of von Neumann algebras in every model of ZF. For example, the easy direction of von Neumann's double commutant theorem implies that the double commutant of any self-adjoint subset of $\cB( \calH)$ is closed in the weak operator topology. More interesting examples of von Neumann algebras are group algebras  (\cite[\S III.3.3]{BlackadarOperator}) and it would be interesting to see how their rich structure behaves in a choiceless setting.  %If in a group~$\Gamma$ all nontrivial conjugacy classes are infinite then $L(\Gamma)$ is a II$_1$ factor (\cite[Proposition~III.3.3.7]{BlackadarOperator}).

An interesting source of Hilbert spaces  is given by the following.

\paragraph{} \begin{Problem}
		For a Boolean algebra $\cB$ with finitely additive probability measure $\mu$ define and study the Hilbert space $L^2(\cB,\mu)$.
	\end{Problem}

The methods of this paper clearly cannot give an example of a Hilbert space that is not included in a Hilbert space with a basis  (see \S\ref{S.universality}), but it is not obvious whether $L^2(\cB,\mu)$ is included in a Hilbert space with a basis for every measure algebra $(\cB,\mu)$.
The results of \cite[\S 56]{fremlin:MT5} may be relevant, and \cite[Exercise~561Y (c) and (i)]{fremlin:MT5} gives a glimpse into the richness of spaces of this sort. See \cite{farah2023choice} for an elaboration on this and reformulations of the assertion that the union of a countable collection of finite sets is countable.

\paragraph{}
\begin{Problem}
Compute the $K$-theory and nonstable $K$-theory of $\cB( \calH)$ and $\cQ( \calH)$ for various
Hilbert spaces $ \calH$.
\end{Problem}

 The answer will be dependent on the model of ZF used.

\bibliography{choiceless-hilbert}
\bibliographystyle{alpha}
\end{document}